\documentclass[11pt]{article}\textwidth 160mm\textheight 235mm
\usepackage{amsfonts}
\usepackage{amssymb}
\usepackage{graphicx}
\usepackage{cancel}
\usepackage{tikz}
\usetikzlibrary{matrix}
\usetikzlibrary{arrows,shapes}
\usepackage{pgfplots}
\pgfplotsset{my style/.append style={axis x line=middle, axis y line=
middle, axis equal }}

\def\ve{\varepsilon}

\def\hat{\widehat}
\def\tilde{\widetilde}
\def\emp{\emptyset}

\def\dom{{\rm dom}\,}
\def\span{{\rm span}\,}
\def\epi{{\rm epi\,}}

\def\lra{\longrightarrow}

\def\N{{\cal N}}

\def\O{{\cal O}}

\def\sub{\partial}
\def\B{\mathbb B}

\def\ox{\overline{x}}
\def\oy{\overline{y}}

\def\cl{{\rm cl}\,}

\def\disp{\displaystyle}

\def\Limsup{\mathop{{\rm Lim}\,{\rm sup}}}

\def\tto{\;{\lower 1pt\hbox{$\rightarrow$}}\kern-10pt
\hbox{\raise 2pt\hbox{$\rightarrow$}}\;}
\def\Hat{\widehat}

\def\Bar{\overline}
\def\ra{\rangle}
\def\la{\langle}
\def\ve{\varepsilon}
\def\B{I\!\!B}
\def\h{\hfill\Box}
\def\R{\mathbb{R}}
\def\N{I\!\!N}
\def\ox{\bar{x}}
\def\oy{\bar{y}}

\def\ov{\bar{v}}

\def\ou{\bar{u}}

\def\co{\mbox{\rm co}}
\def\cone{\mbox{\rm cone}}

\def\int{\mbox{\rm int}\,}
\def\gph{\mbox{\rm gph}\,}
\def\epi{\mbox{\rm epi}\,}

\def\dom{\mbox{\rm dom}\,}
\def\bd{\mbox{\rm bd}\,}

\def\sgn{\mbox{\rm sgn}\,}

\def\cl{\mbox{\rm cl}\,}

\def\h{\hfill\triangle}

\def\O{\Omega}
\def\ph{\varphi}
\def\emp{\emptyset}
\def\st{\stackrel}
\def\oR{\Bar{\R}}
\def\lm{\lambda}
\def\gg{\gamma}
\def\dd{\delta}
\def\al{\alpha}

\def \N{I\!\!N}
\def\th{\theta}

\def\sce{\setcounter{equation}{0}}
\begin{document}
\vspace*{0.5in}
\begin{center}
{\bf GENERALIZED DIFFERENTIATION OF PIECEWISE LINEAR FUNCTIONS IN SECOND-ORDER VARIATIONAL ANALYSIS}\footnote{This research was partly supported by the National Science Foundation under grants DMS-1007132 and DMS-1512846 and by the Air Force Office of Scientific Research grant \#15RT0462.}\\[3ex]
BORIS S. MORDUKHOVICH\footnote{Department of Mathematics, Wayne State University, Detroit, Michigan 48202, USA
(boris@math.wayne.edu).} and M. EBRAHIM SARABI\footnote{Department of Mathematics, Wayne State University, Detroit, Michigan 48202, USA (ebrahim.sarabi@wayne.edu).}
\end{center}
\vspace*{0.05in}

\small{\bf Abstract.} The paper is devoted to a comprehensive second-order study of a remarkable class of convex extended-real-valued functions that is highly important in many aspects of nonlinear and variational analysis, specifically those related to optimization and stability. This class consists of lower semicontinuous functions with possibly infinite values on finite-dimensional spaces, which are labeled as ``piecewise linear" ones and can be equivalently described via the convexity of their epigraphs. In this the paper we calculate the second-order subdifferentials (generalized Hessians) of arbitrary convex piecewise linear functions, together with the corresponding geometric objects, entirely in terms of their initial data. The obtained formulas allow us, in particular, to justify a new exact (equality-type) second-order sum rule for such functions in the general nonsmooth setting.



\newtheorem{Theorem}{Theorem}[section]
\newtheorem{Proposition}[Theorem]{Proposition}
\newtheorem{Remark}[Theorem]{Remark}
\newtheorem{Lemma}[Theorem]{Lemma}
\newtheorem{Corollary}[Theorem]{Corollary}
\newtheorem{Definition}[Theorem]{Definition}
\newtheorem{Example}[Theorem]{Example}
\renewcommand{\theequation}{{\thesection}.\arabic{equation}}
\renewcommand{\thefootnote}{\fnsymbol{footnote}}

\normalsize
\section{Introduction}\sce

Variational analysis has been well recognized as a rapidly developed area of nonlinear analysis, which particularly addresses optimization-related and equilibrium problems, control systems governed by ODEs and PDEs as well as their numerous applications. Furthermore, this broad area of research also deals with many nonvariational issues (stability and generalized differential calculus are among them), which can be resolved by employing variational principles and techniques; see, e.g., the books \cite{bz,m06,rw} and the references therein.

Since nonsmooth functions, sets with nonsmooth boundaries, and set-valued mappings (multifunctions) naturally and frequently appear in variational analysis even for problems with nonsmooth initial data, concepts and machinery of generalized differentiation play a fundamental role in many aspects of variational theory and applications. Extended-real-valued functions, which are usually assumed to be lower semicontinuous (l.s.c.), belong to the central objects of variational analysis giving us a convenient framework for the unified study of functions and sets, combining analytic and geometric methods of their investigation, and having various applications to constrained optimization, stability, and other variational and nonvariational issues.

In this paper we pay the main attention to the study of convex and l.s.c.\ extended-real-valued functions $\ph\colon\R^n\to\oR:=(-\infty,\infty]$, which can be defined geometrically as those whose epigraphs $\epi\ph:=\{(x,\mu)\in\R^{n+1}|\;\mu\ge\ph(x)\}$ are convex polyhedra/polyhedral sets; see Section~3 for analytic descriptions and additional properties. This class was introduced by Rockafellar \cite{r88}, even without the convexity requirement, under the name of {\em piecewise linear} functions in the framework of the so-called epi-subdifferentiability of extended-real-valued functions. Over the years, it has been recognized the importance of piecewise linear functions for their own sake and also as a key building block in the construction of {\em fully amenable} functions playing a crucial role in second-order variational analysis and its numerous applications to constrained optimization and related topics; see \cite{rw} and the subsequent publications listed, e.g., in \cite{mr}.

Our major focus here is on studying the {\em second-order subdifferentials} (or generalized Hessians) of  convex {piecewise} linear functions in the sense initiated by Mordukhovich \cite{m92} for general extended-real-valued functions; see Section~\ref{tools} for the precise definitions. These constructions and their modification have been widely used in second-order variational analysis and its applications to stability issues and necessary optimality conditions for various classes of optimization-related and optimal control problems as well as for parametric variational and equilibrium systems. We refer the reader to \cite{m92,m94,m06,mo01,pr} for the original motivations and underlying results involving the second-order subdifferentials and also to \cite{ch,chhm,dsy,dl,eh,hmn,hos,hr,hy,mnn,mn,mr,ms14,or,os,yy} and the bibliographies therein for more recent studies and applications.

The effective implementation of the obtained general results requires {\em explicit calculations} of the second-order subdifferentials entirely via the problem data. Such calculations can be found in many publications, including the ones listed above. Among those given for particular subcollections of convex piecewise linear functions, we mention calculating the second-order subdifferentials for the indicator functions of convex polyhedra presented in \cite{dr96} via the explicit but not easily implementable ``critical face" condition and in \cite{hmn,hr,yy} obtained in significantly more effective terms involving index sets associated with polyhedral constraints. Quite recently, the complete second-order calculations have been done in \cite{e14,eh} for certain subclasses of the maximum functions, which also belong to the collection of all (convex) piecewise linear ones.

Of course, the calculations obtained for particular classes of functions can be extended to their various combinations provided the availability of appropriate {\em second-order calculus rules}. We distinguish between the {\em exact} (equality-type) calculus rules giving us precise expressions of the second-order subdifferentials for the function combinations under consideration via those of the components involved vs.\ the corresponding inclusions ensuring only {\em upper estimates}. While the second-order subdifferential sum and chain rules of both types can be found in \cite{m02,m06,mnn,mo01,mr,or}, we emphasize the following two issues closely related to our subsequent study: {\bf (i)} The exact second-order sum rules are available only when one of the summands is {\em ${\cal C}^2$-smooth} (or a little less: ${\cal C}^1$-smooth with the strictly differentiable derivative at the reference point); see \cite[Proposition~1.121]{m06}. {\bf (ii)} We have the exact second-order chain rule for a subclass of fully amenable compositions with convex {\em piecewise linear} extended-real-valued {\em outer} functions obtained in \cite[Theorem~4.3]{mr} under a certain {\em second-order qualification condition} (SOCQ) involving the second-order subdifferential of the outer function at zero.

The {\em main contribution} of this paper is a precise calculation of the basic second-order subdifferential for the {\em general class} of convex and piecewise linear functions in finite dimensions entirely in terms of their given data. As a consequence of these calculations, we justify the exact {\em second-order sum rule} for such functions, which is the first result of this type allowing all the summands to be {\em nonsmooth} and even extended-real-valued. Furthermore, the given calculations make it possible to express the aforementioned {\em second-order chain rule} for fully amenable composition and the corresponding SOCQ explicitly via the parameters of piecewise linearity. Various specifications of the general results obtained are also discussed. As a by-product of these developments, we calculate some second-order constructions and their geometric associates, which are preliminary in our consideration while playing a major role prior to the performing crucial limiting procedures and being certainly of their own interest.\vspace*{0.02in}

The rest of the paper is organized as follows. Section~\ref{tools} recalls some definitions and facts from generalized differential theory of variational analysis needed for the formulations and proofs of the subsequent results. In Section~\ref{p-lin} we define the class of (convex) piecewise functions, present their equivalent descriptions, and derive their first-order properties needed in what follows.

In Section~\ref{pre} we start a local second-order analysis of piecewise linear functions calculating geometric second-order constructions generated by certain ``prenormals" to the (first-order) subdifferential and its graph, which may be a nonconvex set despite the convexity of piecewise linear functions under consideration. Although the constructions used in this section are nonrobust and do not possess desired calculus rules, the obtained second-order calculations provide useful tools of analysis and eventually lead us to complete calculations of the basic second-order subdifferential by performing appropriate limiting procedures in the subsequent sections.

Section~\ref{2nd-sec} is the culmination of the paper. It presents explicit and easily implementable calculations of the second-order subdifferential  and related constructions for the general class of convex piecewise linear functions entirely via their given parameters. As a product of these calculations, we justify the aforementioned second-order subdifferential sum rule involving such functions. Section~6 is devoted to some specifications and improvements of our general results for the case of the maximum functions while extending in this way the recent results from \cite{e14,eh}.\vspace*{0.01in}

Our notation and terminology are standard in variational analysis and generalized differentiation; see, e.g., \cite{m06,rw}. Recall that, given an nonempty set $\O\subset\R^n$, the notation $\co\O$, $\cone\O$, $\bd\O$, $\int\O$, $\cl\O$, and $\span \O$ stands for the convex hull, conic hull, boundary, interior, and closure of $\O$, as well as the smallest linear subspace containing $\O$, respectively. Furthermore, $x\st{\O}{\to}\ox$ indicates that $x\to\ox$ with $x\in\O$. For a set-valued mapping $F\colon\R^n\tto\R^m$ the symbol
\begin{eqnarray}\label{1.1}
\disp\Limsup_{x\to\ox}F(x):=\Big\{y\in\R^m\Big|\;\exists\,x_k\to\ox,\;y_k\to y\;\mbox{ with }\;y_k\in F(x_k),\;k\in\N:=\{1,2,\ldots\}\Big\}
\end{eqnarray}
signifies the Painlev\'e-Kuratowski outer limit of $F$ as $x\to\ox$.

\section{Tools of Generalized Differentiation}\label{tools}\sce

Based mainly on \cite{m06,rw}, we introduce and briefly discuss here some notions of generalized differentiation in variational analysis widely used in what follows. Starting with geometric constructions, for a nonempty set $\O\subset\R^n$ define the {\em prenormal cone} (known also as the regular or Fr\'echet normal cone) to $\O$ at $x\in\O$ by
\begin{equation}\label{2.1}
\Hat N(x;\O):=\disp\Big\{v\in\R^n\Big|\;\limsup_{u\st{\O}{\to}x}\frac{\la v,u-x\ra}{\|u-x\|}\le 0\Big\}.
\end{equation}
We choose the prefix ``pre" for (\ref{2.1}) and the related analytic constructions to emphasize that, e.g., (\ref{2.1}) is not a proper normal cone, since it is often trivial (reduces to $\{0\}$) at boundary points of simple nonconvex subsets of $\R^2$, which does not correspond to the meaning of normals. The things change when we employ the limiting procedure via (\ref{1.1}) giving us us the {\em normal cone} (known as the limiting, basic, or Mordukhovich one)
\begin{equation}\label{2.2}
N(\ox;\O)=\disp\Limsup_{x\st{\O}{\to}\ox}\Hat N(x;\O)
\end{equation}
to $\O$ at $\ox\in\O$ nontrivial at boundary points. In spite of (in fact due to) its nonconvexity, the normal cone (\ref{2.2}) together with the associated coderivative and subdifferential constructions (see below) enjoy {\em full calculus} induced by variational/extremal principles of variational analysis. Note that the prenormal cone (\ref{2.1}) is convex being dual/polar
\begin{eqnarray}\label{dua}
\Hat N(x;\O)=T(x;\O)^*:=\Big\{v\in\R^n\Big|\;\la v,w\ra\le 0\;\mbox{ for all }\;w\in T(x;\O)\Big\}
\end{eqnarray}
to the (Bouligand-Severi) {\em tangent/contingent cone} $T(x;\O)$ to $\O$ at $x\in\O$ is defined by
\begin{eqnarray}\label{2.5}
T(x;\O):=\Big\{w\in\R^n\Big|\;\exists\,x_k\st{\O}{\to}x,\;\al_k\ge 0\;\mbox{ with }\;\al_k(x_k-x)\to w\;\mbox{ as }\;k\to\infty\Big\},
\end{eqnarray}
while the normal cone (\ref{2.2}) cannot be tangentially generated due to its intrinsic nonconvexity. If $\O$ is convex, both cones (\ref{2.1}) and (\ref{2.2}) reduce to the classical normal cone of convex analysis.

Given an extended-real-valued function $\ph\colon\R^n\to\oR$ with $\dom\ph:=\{x\in\R^n|\;\ph(x)<\infty\}$, we define the (first-order) {\em subdifferential} of $\ph$ at $\ox\in\dom\ph$ by
\begin{equation}\label{2.6}
\partial\ph(\ox):=\Big\{v\in\R^n\Big|\;(v,-1)\in N\big((\ox,\ph(\ox));\epi\ph\big)\Big\}
\end{equation}
via the normal cone (\ref{2.2}) to its epigraph; see \cite{m06,rw} for the equivalent analytic descriptions.
Note that $\partial\dd(\ox;\O)=N(\ox;\O)$, $\ox\in\O$, for the indicator function $\dd(\cdot;\O)=\dd_\O(\cdot)$ of $\O$ defined by $\dd(x;\O):=0$ if $x\in\O$ and $\dd(x;\O):=\infty$ if $x\notin\O$.

Considering next an arbitrary set-valued mapping $F\colon\R^n\tto\R^m$ with the domain and graph
$$
\dom F\colon=\Big\{x\in\R^n\Big|\;F(x)\ne\emp\Big\},\quad\gph F\colon=\Big\{(x,y)\in\R^n\times\R^m\Big|\;x\in F(x)\Big\},
$$
define its {\em precoderivative} and {\em coderivative} at $(\ox,\oy)\in\gph F$ by respectively,
\begin{equation}\label{2.7}
\Hat D^*F(\ox,\oy)(u):=\Big\{v\in\R^n\Big|\;(v,-u)\in\Hat N\big((\ox,\oy);\gph F\big)\Big\},\quad u\in\R^m,
\end{equation}
\begin{equation}\label{2.8}
D^*F(\ox,\oy)(u):=\Big\{v\in\R^n\Big|\;(v,-u)\in N\big((\ox,\oy);\gph F\big)\Big\},\quad u\in\R^m.
\end{equation}
If $F=f\colon\R^n\to\R^m$ is single-valued and ${\cal C}^1$-smooth around $\ox$, then we have
$$
\Hat D^*f(\ox)(u)=D^*f(\ox)(u)=\Big\{\nabla f(\ox)^*u\Big\}\;\mbox{ for all }\;u\in\R^m,
$$
where $\oy=f(\ox)$ is omitted and the sign $^*$ signifies the matrix transposition.

Following now the dual-space ``derivative-of-derivative" approach to second-order generalized differentiation \cite{m92} of extended-real-valued functions $\ph\colon\R^n\to\oR$, we define the {\em second-order subdifferential} of $\ph$ at $\ox\in\dom\ph$ relative to $\ov\in\partial\ph(\ox)$ by
\begin{equation}\label{secseb}
\sub^2\ph(\ox,\ov)(u)\colon=(D^*\sub\ph)(\ox,\ov)(u),\quad u\in\R^n,
\end{equation}
which is the main object of our study. If $\ph$ is ${\cal C}^2$-smooth around $\ox$, then
$$
\partial^2\ph(\ox)(u)=\Big\{\nabla^2\ph(\ox)u\Big\}\;\mbox{ for all }\;u\in\R^n
$$
via the (symmetric) Hessian $\nabla^2\ph(\ox)$, and so (\ref{secseb}) is viewed as a generalized Hessian of $\ph$ at $\ox$.

\section{Piecewise Linear Functions via First-Order Study}\label{p-lin}\sce

Recall \cite{r88,rw} that $\ph\colon\R^n\to\oR$ is {\em piecewise linear} if its domain $\dom\ph$ is nonempty and can be represented as the union of finitely many convex polyhedral sets so that on each of these pieces $\ph$ is given by $\la a,x\ra-\alpha$ with some $\alpha\in\R$ and $a\in\R^n$.
Observing that such functions are not necessarily convex, we focus on the study of {\em convex} piecewise linear (CPWL) functions, which admit the following equivalent descriptions \cite[Theorem~2.49]{rw}.

\begin{Proposition}{\bf(convex piecewise linear functions).}\label{cpwl} The following are equivalent:

{\bf (i)} $\ph\colon\R^n\to\oR$ is a convex and piecewise linear function labeled as $\ph\in CPWL$.

{\bf (ii)} The epigraph $\epi\ph$ is a convex polyhedron in $\R^{n+1}$.

{\bf (iii)} There are $\alpha_i\in\R$ and $a_i\in\R^n$ for $i\in T_1\colon=\{1,\ldots,l\}$ such that $\ph$ is represented by
\begin{equation}\label{eq00}
\ph(x)=\left\{\begin{array}{ll}
\max\Big\{\la a_1,x\ra-\alpha_1,\ldots,\la a_l,x\ra-\alpha_l\Big\}&\mbox{if }\;x\in\dom\ph,\\
\infty&\mbox{otherwise}
\end{array}\right.
\end{equation}
with some $l\in\N$, where the set $\dom\ph$ is a convex polyhedron given by
\begin{equation}\label{dom}
\dom\ph=\Big\{x\in\R^n\Big|\;\la d_i,x\ra\le\beta_i\;\mbox{ for all }\;i\in T_2:=\{1,\ldots,m\}\Big\}
\end{equation}
with some $d_i\in\R^n$, $\beta_i\in\R$, and $m\in\N$.
\end{Proposition}

It follows from (\ref{eq00}) that any $\ph\in CPWL$ can be represented in the {\em sum form}
\begin{equation}\label{theta}
\ph(x)=\max\Big\{\la a_1,x\ra-\alpha_1,\ldots,\la a_l,x\ra-\alpha_l\Big\}+\dd(x;\dom\ph),\quad x\in\R^n,
\end{equation}
where both summands are nonsmooth. Note that CPWL functions may be given in other forms different from (\ref{theta}), e.g., as the support function of a convex polyhedron
$$
\ph(x)=\sigma_{P}(x):=\sup\{\la p,x\ra|\;p\in P\},
$$
which is conjugate to the indicator function of $P$. Thus $\sigma_P$ is CPWL by \cite[Theorem~11.14(a)]{rw}.\vspace*{0.03in}

The next simple while important observation, giving in particular some relationships between the parameters in (\ref{eq00}) and (\ref{dom}), easily follows from Proposition~\ref{cpwl}.

\begin{Proposition}{\bf (domain of CPWL functions).}\label{srpw} Let $\ph$ be CPWL. Then we have:

{\bf (i)} $\dom\ph=\bigcup^{l}_{i=1} {C_i}$, where $l$ is taken from {\rm(\ref{eq00})} and the sets $C_i$, $i\in T_i$, are defined by
\begin{equation}\label{pwlr1}
 C_i:=\Big\{x\in\dom\ph\Big|\;\la a_j,x\ra-\al_j\le\la a_i,x\ra-\al_i\;\;\mbox{for all}\;\;j\in T_1\Big\}.
\end{equation}

{\bf (ii)} $\ph$ can be equivalently represented via the pieces of linear functions on $C_i$:
\begin{eqnarray*}
\ph(x)=\la a_i,x\ra-\alpha_i\;\mbox{ for all }\;x\in C_i,\;i\in T_1.
\end{eqnarray*}
\end{Proposition}
{\bf Proof}. The inclusion $\bigcup^{l}_{i=1}{C_i}\subset\dom\ph$ follows from the construction of $C_i$. The reverse inclusion is a consequence of (\ref{eq00}). The assertion in (ii) comes directly from (\ref{pwlr1}). $\h$\vspace*{0.05in}

To proceed further, take $\ph\in CPWL$ with $\ox\in\dom\ph$ and consider the index set
\begin{equation}\label{active2}
K(\ox)\colon=\Big\{i\in T_1\Big|\;\ox\in C_i\Big\},
\end{equation}
where each $C_i$ is defined in (\ref{pwlr1}). It is easy to see that
\begin{equation}\label{eq03}
\la a_j,\ox\ra-\al_j=\la a_i,\ox\ra-\al_i\;\mbox{ for any }\;i,j\in K(\ox).
\end{equation}

The next proposition collects basic first-order subgradient facts for CPWL functions.

\begin{Proposition}{\bf(first-order subdifferential of convex piecewise linear functions).}\label{fop} Let $\ph\in CPWL$ with $\ox\in\dom\ph$, and let $K(\ox)$ be defined in {\rm(\ref{active2})}. The following hold:

{\bf (i)}  For any $x$ sufficiently close to $\ox$ we have $\partial\ph(x)\subset\partial\ph(\ox)$.

{\bf (ii)} The first-order subdifferential of $\ph$ at $\ox$ is represented by
\begin{equation}\label{fos}
\partial\ph(\ox)=\co\Big\{a_i\;\Big|\;i\in K(\ox)\Big\}+N(\ox;\dom\ph).
\end{equation}
\end{Proposition}
{\bf Proof}. The subdifferential inclusion in (i) follows directly from the piecewise linear representation of Proposition~\ref{srpw}(ii) with $C_i$ taken from (\ref{pwlr1}). To justify (\ref{fos}), apply to (\ref{theta}) the classical subdifferential sum rule of convex analysis and then the subdifferential formula for the maximum function in (\ref{eq00}), where (\ref{active2}) is the set of the corresponding active indexes. $\h$\vspace*{0.05in}

Note that the normal cone in (\ref{fos}) can be calculated explicitly via the given data of the polyhedral set $\dom\ph$ from (\ref{dom}). To do it, consider the set of {\em active domain indexes}
\begin{equation}\label{active}
I(\ox)\colon=\Big\{i\in T_2\Big|\;\la d_i,\ox\ra=\beta_i\Big\}
\end{equation}
for (\ref{dom}) at $\ox\in\dom\ph$ and recall the well-known tangent cone representation
\begin{equation}\label{tanc}
T(\ox;\dom \ph)=\Big\{x\in\R^n\Big|\;\la d_i,\ox\ra\le 0\;\;\mbox{for all}\;\;i\in I(\ox)\Big\},
\end{equation}
which gives in duality the normal cone one to the polyhedral domain
\begin{equation}\label{norc}
N(\ox;\dom\ph)=\Big\{\sum_{i\in I(\ox)}\mu_i d_i\Big|\;\mu_i\ge 0\Big\}.
\end{equation}
Take now $(\ox,\ov)\in\gph\partial\ph$ and deduce from Proposition~\ref{fop}(ii) the representation
\begin{equation}\label{eq04}
\ov=\ov_1+\ov_2\;\mbox{ with some }\;\ov_1\in\co\Big\{a_i\Big|\;i\in K(\ox)\Big\}\;\mbox{ and }\;\ov_2\in N(\ox;\dom\ph),
\end{equation}
which allows us to find by (\ref{norc}) such $\bar\lm_i\ge 0$ and $\bar\mu_i\ge 0$ that
\begin{equation}\label{eq06}
\disp{\ov_1=\sum_{i\in K(\ox)}\bar\lm_ia_i}\;\mbox{ and }\;\disp{\ov_2=\sum_{i\in I(\ox)}\bar\mu_id_i}\;\mbox{ with }\;\disp{\sum_{i\in K(\ox)}{\bar\lm_i}=1}.
\end{equation}
Corresponding to (\ref{eq06}), define the {\em index sets of positive multipliers} by
\begin{equation}\label{eq05}
J_+(\ox,\ov_1):=\Big\{i\in K(\ox)\Big|\;\bar\lm_i>0\Big\}\;\mbox{ and }\;J_+(\ox,\ov_2):=\Big\{i\in I(\ox)\Big|\;\bar\mu_i>0\Big\}.
\end{equation}

The next result gives us an effective description of the points in the subdifferential graph for CPWL functions. It is certainly of its own interest while playing an important role in the subsequent calculations of second-order constructions for such functions. To provide a better understanding of it, consider the diagrams presented in Fig.\ 1 and Fig.\ 2.
\begin{figure}[b]
\begin{minipage}[t]{0.5\textwidth}
\begin{tikzpicture}[yscale=1,xscale=1]
\draw[thick,->](-3,0)--(3,0);
\draw[thick,->](0,-1.5)--(0,3);
\draw[ultra thick,domain=0:1] plot(\x,{(0.5*\x)});
\draw[ultra thick,domain=-2:0] plot(\x,{(-1*\x)});
\draw[ultra thick,domain=1:2] plot(\x,{(2*\x-1.5)});
\draw[dashed](0,.5)node[left]{\scriptsize $0.5$}--(1,.5)
--(1,0)node[below]{\scriptsize $1$};
\draw[dashed](0,2) node[right]{\scriptsize $2$}--(-2,2)
--(-2,0)node[below]{\scriptsize $-2$};
\draw[dashed](0,2.5)node[left]{\scriptsize $2.5$}--(2,2.5)
--(2,0)node[below]{\scriptsize 2};
\draw[<->](0,0)--node[below]{\tiny$ C_1$}(1,0);
\draw[<->](1,0)--node[below]{\tiny$C_2$}(2,0);
\draw[<->](-2,0)--node[below]{\tiny$C_3$}(0,0)node[below=2cm]{Fig.~1: Graph of $\ph$};
\end{tikzpicture}
\end{minipage}
\begin{minipage}[t]{0.5\textwidth}
\begin{tikzpicture}[yscale=1,xscale=1]
\draw[thick,->](-3,0)--(3,0);
\draw[thick,->](0,-1.5)--(0,3);
\draw[ultra thick,domain=1:4]plot(\x,{(2)});
\draw[ultra thick,domain=0:1]plot(\x,{(.5)});
\draw[ultra thick](1,0.5)--(1,2);
\draw[ultra thick](0,0.5)--(0,-1);
\draw[ultra thick,domain=-3:0]plot(\x,{(-1)});
\draw[dashed](0,2)node[left]{\scriptsize$2$}--(1,2);
\draw[dashed](1,.5)--(1,0)node[below]{\scriptsize 1};
\draw[fill](1,1)circle(2pt);
\draw[dashed](0,1)node[left]{\scriptsize$1$}--(1,1)node[right]{\scriptsize $A$};
\draw[dashed](0,.5)node[left]{\scriptsize$0.5$};
\draw[dashed](0,-1)node[right]{\scriptsize$-1$};
\draw[dashed](2,0)node[below]{\scriptsize$2$}--(2,2);
\draw[dashed](-2,0)node[above]{\scriptsize$-2$}--(-2,-1);
\draw[<->](0,0)--node[below]{\tiny$ C_1$}(1,0);
\draw[<->](1,0)--node[below]{\tiny$C_2$}(2,0);
\draw[<->](-2,0)--node[below]{\tiny$C_3$}(0,0)node[below=2cm]
{Fig.\ 2: Graph of $\sub\ph$};
\end{tikzpicture}
\end{minipage}
\end{figure}
Given $A\colon=(\ox,\ov)=(1,1)\in\gph\sub\ph$, we have $\ov=1\in\sub\ph(1)=[\frac{1}{2},2]$ and so $\ov=\bar\lm_1\frac{1}{2}+\bar\lm_2 2$ with $\bar\lm=(\bar\lm_1,\bar\lm_2)=(\frac{2}{3},\frac{1}{3})$. This tells us that $J_+(\ox,\ov_1) =\{1,2\}$, where $\ov=\ov_1$. Suppose that $(x,v)$ is a point of $\gph\sub\ph$, which is sufficiently close to $A$. As can be observed from the graph of $\sub\ph$, the only possibility for converging to $A$ while remaining in the graph of $\sub\ph$ is $(x,v)=(1,1\pm\epsilon)$ when $\epsilon\downarrow 0$. Thus it confirms that $x\in C_1\cap C_2=\{1\}$ and hence $x\in C_i$ for $i\in J_+(\ox,\ov_1)$.\vspace*{0.03in}

Now we are ready to justify this phenomenon for the general class of CPWL functions.

\begin{Theorem}{\bf(description of points in the subdifferential graph of CPWL functions).}\label{gphpar}
Let $\ph\in CPWL$ with $(\ox,\ov)\in\gph\partial\ph$. Then there exists a neighborhood $O$ of $(\ox,\ov)$ such that for any $(x,v)\in(\gph\partial\ph)\cap O$ we have $x\in\bigcap_{i\in J_+(\ox,\ov_1)}{C_i}$, where $\ov_1$ and $J_+(\ox,\ov_1)$ are taken from {\rm(\ref{eq06})} and {\rm(\ref{eq05})}, respectively, and where the polyhedral sets $C_i$ are defined in {\rm (\ref{pwlr1})}.
\end{Theorem}
{\bf Proof}. We split the proof into the following major steps with keeping all the notation above.\vspace*{0.02in}

{\bf Claim~1:} {\em Let $\ov=\sum_{i\in P}\eta_ia_i+\sum_{i\in Q}\tau_id_i$ with some $\tau_i,\eta_i\ge 0$ satisfying $\sum_{i\in P}{\eta_i}=1$, $P\subset K(\ox)$, and $Q\subset I(\ox)$. Then we have the equality
\begin{equation}\label{eq07}
\disp\sum_{i\in P}\eta_i\al_i+\sum_{i\in Q}\tau_i\beta_i=\sum_{i\in K(\ox)}\bar\lm_i\al_i+\sum_{i\in I(\ox)}\bar\mu_i\beta_i,
\end{equation}
where the multipliers $\bar\lm_i$ and $\bar\mu_i$ are taken from {\rm(\ref{eq06})}.}\\[1ex]
To verify this claim, suppose that $\ov=\hat v_1+\hat v_2$ for $\hat v_1=\sum_{i\in P}\eta_i a_i$ and $\hat v_2=\sum_{i\in Q}\tau_id_i$  with $\sum_{i\in P}{\eta_i}=1$ and $\eta_i,\tau_i\ge 0$. Fix $j\in P$ and observe that $\la a_j,\ox\ra-\al_j=\la a_i,\ox\ra-\al_i$ for any $i\in K(\ox)$ due to (\ref{eq03}) and $P\subset K(\ox)$. It tells us that
$$
\la a_j,\ox\ra-\al_j=\la\ov_1,\ox\ra-\sum_{i\in K(\ox)}\bar\lm_i\al_i
$$
with $\ov_1$ taken from (\ref{eq06}), which implies in turn that
\begin{equation}\label{eq08}
\la\hat v_1,\ox\ra-\sum_{j\in P}\eta_j\al_j=\la\ov_1,\ox\ra-\sum_{i\in K(\ox)}\bar\lm_i\al_i.
\end{equation}
Since $Q\subset I(\ox)$, this allows us to deduce that
\begin{equation}\label{eq09}
\la\hat v_2,\ox\ra=\sum_{i\in Q}\tau_i\beta_i\;\mbox{ and }\;\la\ov_2,\ox\ra=\sum_{i\in I(\ox)}\bar\mu_i\beta_i,
\end{equation}
where $\ov_2$ is from (\ref{eq06}). Combining (\ref{eq08}) and (\ref{eq09}) with $\hat v_1+\hat v_2=\ov_1+\ov_2$ justifies the claim.\vspace*{0.03in}

Suppose now that the conclusion of the theorem does not hold and thus find a sequence $(x_k,v_k)\in\gph\sub\ph$ such that
$(x_k,v_k)\to(\ox,\ov)$ as $k\to\infty$ while $x_k\not\in\bigcap_{i\in J_+(\ox,\ov_1)}{C_i}$ for all $k\in\N$. Taking into account that there are only {\em finitely many} convex polyhedral sets $C_i$ and considering a subsequence of $x_k$ if necessary, suppose that $x_k\not\in C_s$ for some $s\in J_+(\ox,\ov_1)$. Furthermore, it is not hard to see that $K(x_k)\subset K(\ox)$ for $k$ sufficiently large. Extracting similarly another subsequence, pick without loss of generality a constant index subset $P\subset K(\ox)$ so that $K(x_k)=P$ for all $k$. Select $j\in P$ and observe that $x_k\in C_j$, which implies by (\ref{pwlr1}) that
\begin{equation}\label{eq010}
\la a_j,x_k\ra-\al_j\ge\la a_i,x_k\ra-\al_i\;\mbox{ for all }\;i\in K(\ox).
\end{equation}
On the other hand, the construction in (\ref{pwlr1}) and the conditions $x_k\not\in C_s$, $x_k\in\dom\ph$ allow us to select $t\in T_1$ independently of $k$ and  so that $\la a_s,x_k\ra-\al_s<\la a_t,x_k\ra-\al_t$.\vspace*{0.02in}

{\bf Claim~2:} {\em We have $t\in K(\ox)$ for the index $t\in T_1$ selected above}.\\[1ex]
Indeed, suppose by contradiction that $t\not\in K(\ox)$. Combining this with $s\in K(\ox)$ tells us that $\la a_t,\ox\ra-\al_t<\la a_s,\ox\ra-\al_s$, and thus
$\la a_t,x_k\ra-\al_t<\la a_s,x_k\ra-\al_s$ for all $k$ sufficiently large. This clearly contradicts the choice of the index $t$ and hence justifies the claim.\vspace*{0.03in}

For any fixed $s\in K(\ox)$ define now the index set
\begin{eqnarray*}
D_s:=\Big\{t\in T_1\Big|\;\la a_s,x_k\ra-\al_s<\la a_t,x_k\ra-\al_t\;\mbox{ for all }\;k\in\N\Big\}.
\end{eqnarray*}
It follows from Claim~2 that $\emp\ne D_s\subset K(\ox)$. We continue with the next assertion.\vspace*{0.02in}

{\bf Claim~3:} {\em $P\subset D_s$, where $P$ was selected so that $K(x_k)=P$ for all $k\in\N$.}\\[1ex]
Assuming the contrary, find $j\in P$ such that $j\not\in D_s $ and pick $t\in D_s$. Employing this gives us
\begin{equation}\label{eq012}
\la a_s,x_k\ra-\al_s<\la a_t,x_k\ra-\al_t\;\mbox{ for all }\;k\in\N.
\end{equation}
Since $j\not\in D_s$, there exists a number $k_0\in\N$ for which we have
\begin{equation}\label{eq013}
\la a_j,x_{k_0}\ra-\al_j\le\la a_s,x_{k_0}\ra-\al_s.
\end{equation}
Combining (\ref{eq012}) for $k=k_0$ together with (\ref{eq013}) leads us to the strict inequality
\begin{eqnarray*}
\la a_j,x_{k_0}\ra-\al_j<\la a_t,x_{k_0}\ra-\al_t,
\end{eqnarray*}
which contradicts (\ref{eq010}) due to $t\in K(\ox)$ and thus verifies the claim.\vspace*{0.03in}

To proceed further, extract another subsequence of $x_k$ and find $Q\subset I(\ox)$ with
\begin{eqnarray*}
\la d_i,x_k\ra=\beta_i\;\mbox{ if }\;i\in Q\;\mbox{ and }\;\la d_i,x_k\ra<\beta_i\;\mbox{ if }\;i\in T_2\setminus Q,
\end{eqnarray*}
which shows that $I(x_k)=Q$ for all $k$. The next claim is as follows.\vspace*{0.02in}

{\bf Claim~4:} {\em We have $\ov\not\in\co\Big\{a_i\Big|\;i\in P\;\Big\}+\cone\Big\{d_i\Big|\;i\in Q\;\Big\}$.}\\[1ex]
To verify the claim, suppose on the contrary that there exist vectors $\hat v_1\in{\rm co}\{a_i|\;i\in P\}$ and $\hat v_2\in{\rm cone}\{d_i|\;i\in Q\}$ such that $\ov=\hat v_1+\hat v_2$. This allows us to find numbers $\tau_i,\eta_i\ge 0$ with $\sum_{i\in P}\eta_i=1$ such that $\hat v_1\colon=\sum_{i\in P}{\eta_i a_i}$ and $\hat v_2\colon=\sum_{i\in Q}{\tau_i d_i}$. Pick $t\in D_s\cap P$, which can be done by Claim~3. It follows from $t\in D_s$ that
\begin{eqnarray*}
\la a_s,x_k\ra-\al_s<\la a_t,x_k\ra-\al_t\;\mbox{ for all }\;k
\end{eqnarray*}
while $t\in P$ results in the inequality
\begin{eqnarray*}
\la a_i,x_k\ra-\al_i\le\la a_t,x_k\ra-\al_t\;\mbox{ whenever }\;i\in K(\ox).
\end{eqnarray*}
Using these two facts together with $s\in J_+(\ox,\ov_1)$ yields
\begin{equation}\label{eq018}
\la \ov_1,x_k\ra-\sum_{i\in K(\ox)}\bar\lm_i\al_i<\la a_t,x_k\ra-\al_t,
\end{equation}
where $\ov_1$ is from (\ref{eq04}) and the multipliers $\bar\lm_i$ are taken from (\ref{eq06}). Remembering that $\la a_i,x_k\ra-\al_i=\la a_j,x_k\ra-\al_j$ for all $i,j\in P$ and taking (\ref{eq018}) into account ensure that
\begin{eqnarray}\label{eq019}
\la \ov_1,x_k\ra-\sum_{i\in K(\ox)}\bar\lm_i\al_i<\la \hat v_1,x_k\ra-\sum_{i\in P}\eta_i\al_i.
\end{eqnarray}
On the other hand, we know that $x_k\in\dom\ph$, which leads us to
\begin{equation}\label{eq020}
\la \ov_2,x_k\ra\le\sum_{i\in I(\ox)}\bar\mu_i\beta_i\;\mbox{ and }\;\la\hat v_2,x_k\ra=\sum_{i\in Q}\tau_i\beta_i.
\end{equation}
Using (\ref{eq019}) and (\ref{eq020}) together with $\ov_1+\ov_2=\hat v_1+\hat v_2$ gives us
\begin{equation}\label{eq021}
\sum_{i\in K(\ox)}\bar\lm_i\al_i+\sum_{i\in I(\ox)}\bar\mu_i\beta_i>\sum_{i\in P}\eta_i\al_i+\sum_{i\in Q}\tau_i\beta_i.
\end{equation}
Appealing finally to Claim~1 along with the inclusions $Q\subset I(\ox)$ and $P\subset K(\ox)$, we arrive at a contradiction with (\ref{eq021}) and hence verify this claim.

Now we are ready to finish the proof of the theorem. Remember that $(x_k,v_k)\st{\tiny\gph\partial\ph}{\lra}(\ox,\ov)$, which yields $v_k\in\partial\ph(x_k)$ for all $k\in\N$. It follows from (\ref{fos}) that $\partial\ph(x_k)=\co\{a_i|\;i\in P\}+\cone\{d_i|\;i\in Q\}$ due to $K(x_k)=P$ and $I(x_k)=Q$. Hence we have $\ov\in\co\{a_i|\;i\in P\}+\cone\{d_i|\;i\in Q\}$ thus contradicting Claim~4 and showing that the assumption made after Claim~1 cannot be correct, while the opposite is the conclusion of the theorem. $\h$

\section{Calculating Prenormals to Subdifferential Graphs}\label{pre}\sce

Here we start our second-order analysis for a general class of CPWL functions following the dual-space derivative-of-derivative approach to second-order generalized differentiation \cite{m92} that focuses on considering certain normals to (first-order) subdifferential graphs. The main result of this section provides a precise calculation of the prenormal/regular normal cone to the graph of the subdifferential mapping for an  arbitrary CPWL function in terms of its given data. As a by-product of this, we derive useful formulas for the classical normal and tangent cones to the (convex) set of subgradients for such functions expressed via their parameters from (\ref{eq00}) and (\ref{dom}). Let us begin with calculations of the latter constructions.

\begin{Lemma}{\bf (normal and tangent cones to subgradient sets of CPWL functions).}\label{polnc} Given $\ph\in CPWL$ and $(\ox,\ov)\in\gph\partial\ph$, take $\ov_1$, $\ov_2$ from {\rm(\ref{eq06})} such that $\ov=\ov_1+\ov_2$. Denote by $K\colon=K(\ox)$, $I\colon=I(\ox)$, $J_{1}:=J_+(\ox,\ov_1)$, and $J_{2}:=J_+(\ox,\ov_2)$ the index sets from {\rm(\ref{active2})}, {\rm(\ref{active})}, and {\rm(\ref{eq05})}, respectively. Then we have the following formulas via the data in {\rm (\ref{eq00})} and {\rm(\ref{dom})}:

{\bf (i)} The normal cone to the subgradient set $\partial\ph(\ox)$ at $\ov$ is calculated by:
\begin{equation}\label{eq060}
\begin{array}{ll}
N(\ov;\partial\ph(\ox))=\Big\{u\in\R^n\Big|&\la a_i-a_j,u\ra=0\;\;\mbox{for}\;\;i,j\in J_1,\\
&\la a_i-a_j,u\ra\le 0\;\;\mbox{for}\;\;(i,j)\in(K\setminus J_1)\times J_1,\\
&\la d_i,u\ra=0\;\;{\mbox{for}}\;\;i\in J_2\;\;\mbox{and}\;\;\la d_i,u\ra\le 0\;\;\mbox{for}\;\;i\in I\setminus J_2\;\Big\}.
\end{array}
\end{equation}

{\bf (ii)} The tangent cone to $\partial\ph(\ox)$ at $\ov$ is expressed as
\begin{equation}\label{eq062}
T(\ov;\partial\ph(\ox))=N(\ov;\partial\ph(\ox))^*=T(\ov_1;A(\ox))+N(\ox;\dom\ph)+\R_{-}\{\ov_2\}
\end{equation}
with $N(\ox;\dom\ph)$ from {\rm(\ref{eq062})}, $A(\ox)\colon=\co\{a_i|\;i\in K(\ox)\}$, and $T(\ov_1;{A}(\ox))$ calculated by
\begin{equation}\label{eq061}
T(\ov_1;A(\ox))=\span\Big\{a_i-a_j\Big|\;i,j\in J_1\Big\}+\cone\Big\{a_i-a_j\Big|\;(i,j)\in(K\setminus J_1)\times J_1\Big\}.
\end{equation}
Finally, we have the explicit formula for calculating the tangent cone to $\partial\ph(\ox)$ at $\ov$:
\begin{equation}\label{eq074}
\begin{array}{lll}
T(\ov;\partial\ph(\ox))&=&\span\Big\{a_i-a_j\Big|\;i,j\in J_1\Big\}+\cone\Big\{a_i-a_j\Big|\;(i,j)\in(K\setminus J_1)\times J_1\Big\}\\
&+&\span\Big\{d_i\Big|\;i\in J_2\Big\}+\cone\Big\{d_i\Big|\;i\in I\setminus J_2\Big\}.
\end{array}
\end{equation}
\end{Lemma}
{\bf Proof}. Proposition~\ref{fop}(ii) allows us to represent $\partial\ph(\ox)={A}(\ox)+{B}(\ox)$ with ${A}(\ox)$ from above and $B(\ox)\colon=N(\ox;\dom\ph)=
\cone\{d_i|\;i\in K(\ox)\}$ by (\ref{norc}). Since ${A}(\ox)$ is convex, we have
\begin{equation}\label{eq065}
\begin{array}{lll}
N(\ov_1;A(\ox))&=&\Big\{u\in\R^n\Big|\;\la u,v_1-\ov_1\ra\le 0\;\;\mbox{for}\;\;v_1\in A(\ox)\Big\},\\
&=&\Big\{u\in\R^n\Big|\;\la u,a_i-a_j\ra=0\;\;\mbox{for}\;\;i,j\in J_1,\\
&&\hspace{1.7cm}\la u,a_i-a_j\ra\le 0\;\;\mbox{for}\;\;(i,j)\in(K\setminus J_1)\times J_1\Big\}.
\end{array}
\end{equation}
Furthermore, it follows from the structure of $B(\ox)$ that
\begin{equation}\label{eq063}
N(\ov_2;B(\ox))=B(\ox)^*\cap\{\ov_2\}^\bot=T(\ox;\dom\ph)\cap\{\ov_2\}^\bot,
\end{equation}
which in turn leads us by (\ref{tanc}) to the representation
\begin{equation}\label{eq064}
N(\ov_2;{B}(\ox))=\Big\{u\Big|\;\la d_i,u\ra= 0\;\;{\mbox{for}}\;\;i\in J_2\;\;\mbox{and}\;\;\la d_i,u\ra \le 0\;\;\mbox{for}\;\;i\in I\setminus J_2\;\Big\}.
\end{equation}
Employing the well-known formula for normals to set additions (see, e.g., \cite[Exercise~6.44]{rw}) shows that $N(\ov;\partial\ph(\ox))=N(\ov_1;A(\ox))\cap N(\ov_2;B(\ox))$. This finally brings us to (\ref{eq060}) due to (\ref{eq065}) and (\ref{eq064}) and thus justifies the normal cone formula (\ref{eq060}) in assertion (i).

To verify now (\ref{eq062}) in (ii), we apply \cite[Corollary~19.3.3]{r70} to the set intersection in (\ref{eq063}) and get by \cite[Exercise~6.44]{rw} and the classical duality between the normal and tangent cones that
\begin{eqnarray*}
\begin{array}{lll}
N(\ov;\partial\ph(\ox))^*&=&T(\ov;\partial\ph(\ox))=\cl[T(\ov_1;A(\ox))+T(\ov_2;B(\ox))]\\
&=&T(\ov_1;A(\ox))+T(\ov_2;B(\ox))=T(\ov_1;A(\ox))+N(\ov_2;B(\ox))^*\\
&=&T(\ov_1;A(\ox))+\Big(T(\ox;\dom\ph)\cap\{\ov_2\}^\bot\Big)^*\\
&=&T(\ov_1;A(\ox))+N(\ox;\dom\ph)+\R_{-}\{\ov_2\},
\end{array}
\end{eqnarray*}
where the last equality follows from \cite[Corollary~16.4.2]{r70} and the choice of $\ov_2\in N(\ox;\dom\ph)$. The tangent cone formula (\ref{eq061}) follows from immediately from (\ref{eq065}) by \cite[Lemma~6.45]{rw}.

To establish finally representation (\ref{eq074}), it remains to use the tangent sum rule $T(\ov;\partial\ph(\ox))=T(\ov_1;A(\ox))+T(\ov_2;B(\ox))$ and calculate $T(\ov_2;B(\ox))$. By the classical Farkas lemma (see, e.g., \cite[Lemma~6.45]{rw}) we deduce from (\ref{eq064}) that
$$
T(\ov_2;B(\ox))=\span\Big\{d_i\Big|\;i\in J_2\Big\}+\cone\Big\{d_i\Big|\; i\in I\setminus J_2\Big\}
$$
and thus complete the proof of the lemma. $\h$\vspace*{0.05in}

The next lemma of its independent interest is significantly more involved providing the representation of the prenormal cone (\ref{2.1}) to the nonconvex subdifferential graphs via the tangent and normal cones to the convex subgradient sets for arbitrary CPWL functions. It reduces to \cite[Proposition~3.2]{hmn} in the case of indicator functions of convex polyhedra, while in the general case the proof is based on the subdifferential graph description established in Theorem~\ref{gphpar}.

\begin{Lemma}{\bf(prenormals to subdifferential graphs via tangents and normals to subgradient sets for CPWL functions).}\label{prenormal}
In the setting of Lemma~{\rm\ref{polnc}} we have the representation
\begin{equation}\label{pncg}
\Hat{N}((\ox,\ov);\gph\partial\ph)=T(\ov;\partial\ph(\ox))\times N(\ov;\partial\ph(\ox)).
\end{equation}
\end{Lemma}
{\bf Proof}. We designate in the proof of the lemma several major steps as follows. Aiming first to justify the inclusion``$\subset$" in (\ref{pncg}), fix $(w,u)\in\Hat{N}((\ox,\ov);\gph\partial\ph)$ telling us by (\ref{2.1}) that
\begin{equation}\label{eq040}
\limsup_{(x,v)\st{\tiny\gph\partial\ph}{\to}(\ox,\ov)}\frac{\la w,x-\ox\ra+\la u,v-\ov\ra}{\|x-\ox\|+\|v-\ov\|}\le 0.
\end{equation}
Letting $x=\ox$ and $v\in\partial\ph(\ox)$ in (\ref{eq040}) tells us that $u\in\Hat N(\ov;\partial\ph(\ox))=N(\ov;\partial\ph(\ox))$. To justify the aforementioned inclusion in (\ref{eq040}), it remains to show that $w\in T(\ov;\partial\ph(\ox))=N(\ov;\partial\ph(\ox))^*$, which amounts to verifying the inequality
\begin{equation}\label{eq041}
\la w,p\ra\le 0\;\mbox{ whenever }\;p\in N(\ov;\partial\ph(\ox)).
\end{equation}
To proceed, pick $p\in N(\ov;\partial\ph(\ox))$, let $x_k:=\ox+\frac{1}{k}p$, $k\in\N$, and then prove the following fact. \vspace*{0.02in}

{\bf Claim~1:} {\em We have $J_2\subset I(x_k)$ and $x_k\in\bigcap_{i\in J_1}{C_i}$, $k\in\N$, with the sets $C_i$ defined in {\rm(\ref{pwlr1})}.}\\[1ex]
To check first that $x_k\in\dom\ph$ as $k\in\N$, we need showing by (\ref{dom}) that $\la d_i,x_k\ra\le\beta_i$ for all $i\in T_2$. Pick the generating vector $d_i$ with $i\in I$ and observe that $\hat v:=\ov+d_i\in\partial\ph(\ox)$ due to (\ref{fos}). Since $p\in N(\ov;\partial\ph(\ox))$, it follows that
\begin{equation}\label{eq043}
\la p,d_i\ra=\la p,\hat v-\ov\ra\le 0.
\end{equation}
Employing (\ref{eq043}) together with $\ox\in\dom\ph$ yields
$$
\la d_i,x_k\ra=\la d_i,\ox\ra+k^{-1}\la d_i,p\ra\le\beta_i+k^{-1}(0)=\beta_i\;\mbox{ whenever }\;i\in I.
$$
The latter inequality holds also for $i\in T_2\setminus I$ by $\la d_i,\ox\ra<\beta_i$ and $x_k\to\ox$. Thus we get $x_k\in\dom\ph$.

Further, pick $t\in J_1$ and show that $x_k\in C_t$, which means that $\la a_j,x_k\ra-\al_j\le\la a_t,x_k\ra-\al_t$, or equivalently $\la a_j-a_t,x_k\ra\le\al_j-\al_t$ for any $j\in T_1$. Indeed, for $j\in T_1\setminus K$ the latter inequality is a consequence of $\la a_j-a_t,\ox\ra<\al_j-\al_t$ and $x_k\to\ox$. Observe otherwise that
\begin{equation}\label{eq042}
\la a_j- a_t ,\ox\ra=\al_j-\al_t\;\mbox{ for }\;j\in K.
\end{equation}
It follows from the representation of $\ov_1$ in (\ref{eq06}) that $\bar\lm_t>0$ therein if $t\in J_1$. Take now $0<\ve <\bar\lm_t$ and get $\tilde v:=\sum_{i\in K}\lm_i a_i\in\partial\ph(\ox)$ with
$$
\lm_i:=\left\{\begin{array}{ll}
\bar\lm_t-\ve&\mbox{if}\;i=t,\\
\bar\lm_j+\ve&\mbox{if}\;i=j,\\
\bar\lm_i&\mbox{otherwise.}
\end{array}\right.
$$
Since $p\in N(\ov;\partial\ph(\ox))$, we get $\la p,\ve(a_j-a_t)\ra=\la p,\tilde v-\ov\ra\le 0$ and therefore arrive at $\la a_j-a_t,p\ra\le 0$. This allows us to obtain the equalities
$$
\la a_j-a_t,x_k-\ox\ra=\la a_j-a_t,k^{-1}p\ra\le 0\;\mbox{ for all }\;j\in K,\;t\in J_1,\;\mbox{ and }\;k\in\N,
$$
which in turn imply that $\la a_j-a_t,x_k\ra\le\la a_j-a_t,\ox\ra=\al_j-\al_t$ by (\ref{eq042}) and thus verifies $x_k\in C_t$.

To justifies the claim, it remains to show that $J_2\subset I(x_k)$, which amounts to saying that $\la p,d_t\ra=0$ for any $t\in J_2$. Noting that the inequality $\la p,d_t\ra\le 0$ follows from (\ref{eq043}) due to $J_2\subset I$, we prove now the converse inequality. Observe that $\breve{v}:=\ov_1+\sum_{i\in I}{\mu_id_i}\in\partial\ph(\ox)$, where the multipliers $\mu_i$ are defined via $\bar\mu_i$ from (\ref{eq06}) by
$$
\mu_i:=\left\{\begin{array}{ll}
\disp\frac{1}{2}\bar\mu_t&\mbox{if }\;i=t,\\
\bar\mu_i&\mbox{otherwise.}
\end{array}\right.
$$
Since $p\in N(\ov;\partial\ph(\ox))$, we get $\la p,-\frac{1}{2}\bar\mu_t d_t\ra=\la p,\breve{v}-\ov\ra\le 0$ yielding $\la p,d_t\ra\ge 0$ by $\bar\mu_t>0$. It shows that $\la p,d_t\ra=0$, and so $t\in I(x_k)$, which finishes the proof of this claim.\vspace*{0.03in}

Note that it follows from Claim~1 by the definitions of $J_{1}=J_+(\ox,\ov_1)$ in (\ref{eq05}) and $K(x)$ in (\ref{active2}) that $J_1\subset K(x_k)$ for all $k$. Since $\ov\in\co\{a_i|\;i\in J_1\}+\cone\{d_i|\;i\in J_2\}$ by the above, we have $\ov\in\partial\ph(x_k)$, $k\in\N$. Recalling that $(x_k,\ov)\st{\tiny \gph\partial\ph}\lra(\ox,\ov)$ and substituting $(x_k,\ov)$ into (\ref{eq040}) for $(x,v)$, we arrive at $\la w,p\ra\le 0$, which justifies (\ref{eq041}) and hence the inclusion ``$\subset$" in (\ref{pncg}).\vspace*{0.03in}

To verify the opposite inclusion in (\ref{pncg}), pick $(w,u)\in T(\ov;\partial\ph(\ox))\times N(\ov;\partial\ph(\ox))$ and show that it satisfies (\ref{eq040}). For any $(\hat x,\hat v)\in\gph\sub\ph$ sufficiently close to $(\ox,\ov)$, observe that $\hat v\in\sub\ph(\ox)$ by Proposition~\ref{fop}(i). Then $u\in N(\ov;\partial\ph(\ox))$ yields $\la u,\hat v-\ov\ra\le 0$ by the convexity of $\partial\ph(\ox)$. Lemma~\ref{polnc}(ii) gives us $w_1\in T(\ov_1;A(\ox))$ with $A(\ox)=\co\{a_i|\;i\in K(\ox)\}$, $w_2\in N(\ox;\dom\ph)$, and $\gg\ge 0$ such that $w=w_1+w_2-\gg\ov_2$. We proceed further with verifying the following fact.\vspace*{0.02in}

{\bf Claim~2:} {\em We have $\la w_1,\hat x-\ox\ra\le 0$ whenever $(\hat x,\hat v)\in\gph\sub\ph$ is sufficiently close to $(\ox,\ov)$}\\[1ex]
To justify this claim, it suffices to show due to (\ref{eq061}) that $\la a_i-a_j,\hat x-\ox\ra=0$ when $i,j\in J_1$ and that $\la a_i-a_j,\hat x-\ox\ra\le 0$
when $(i,j)\in(K\setminus J_1)\times J_1$. In the former case we get from (\ref{eq03}) that
\begin{equation}\label{eq070}
\la a_i-a_j,\ox\ra=\al_i-\al_j\;\mbox{ for all }\;i,j\in J_1.
\end{equation}
Furthermore, $\hat x\in\bigcap_{i\in J_1}{C_i}$ by Theorem~\ref{gphpar}. This shows that $J_1\subset K(\hat x)$ yielding in turn $i,j\in K(\hat x)$. Thus $\la a_i,\hat x\ra-\al_i=\la a_j,\hat x\ra-\al_j$, which together with (\ref{eq070}) leads us to $\la a_i-a_j,\hat x-\ox\ra=0$.

To verify the remaining part of the claim, take  $(i,j)\in(K\setminus J_1)\times J_1$ and deduce from Theorem~\ref{gphpar} that $\hat x\in C_j$. Hence
$\la a_i,\hat x\ra-\al_i\le\la a_j,\hat x\ra-\al_j$ by (\ref{pwlr1}) and therefore we get $\la a_i-a_j,\hat x-\ox\ra\le 0$ by (\ref{eq070}), which completes the proof of the claim.\vspace*{0.03in}

The next claim gives us the final estimate needed in the lemma.\vspace*{0.02in}

{\bf Claim~3:} {\em We have the inequality $\la w_2-\gg\ov_2,\hat x-\ox\ra\le\gg\la\hat v-\ov,\hat x-\ox\ra$.}\\[1ex]
To verify it, observe that $\hat x\in\dom\ph $ and thus get $\la w_2,\hat x-\ox\ra\le 0$ by $w_2\in N(\ox;\dom\ph)$. It follows from the convexity of
$\ph$ that $\sub\ph$ is a monotone mapping. This tells us that $\la\hat v-a_i,\hat x-\ox\ra\ge 0$ due to $\hat v\in\sub\ph(\hat x)$ and $a_i\in\sub\ph(\ox)$ for all $i\in K$. Thus we arrive at
\begin{eqnarray*}
\la a_i-\ov,\hat x-\ox\ra\le\la\hat v-\ov,\hat x-\ox\ra\;\mbox{ for all }\;i\in K,
\end{eqnarray*}
which in turn leads us to the following relationships:
\begin{eqnarray*}
\begin{array}{ll}
\la w_2-\gg\ov_2,\hat x-\ox\ra&\le-\la\gg\ov_2,\hat x-\ox \ra=\gg\la\ov_1-\ov,\hat x-\ox\ra=\gg\disp\sum_{i\in K}\bar\lm_i\la a_i-\ov,\hat x-\ox\ra\\
&\le\disp\gg\sum_{i\in K}\bar\lm_i\la\hat v-\ov,\hat x-\ox\ra=\gg\la\hat v-\ov,\hat x-\ox\ra
\end{array}
\end{eqnarray*}
and thus justifies the claimed inequality.\vspace*{0.03in}

To complete finally the proof of the lemma, we combine Claim~2 and Claim~3 together with the inequality $\la u,\hat v-\ov\ra\le 0$ obtained above. This allows us to conclude that
$$
\frac{\la w,\hat x-\ox\ra+\la u,\hat v-\ov\ra}{\|\hat x-\ox\|+\|\hat v-\ov\|}\le\frac{\gg\|\hat v-\ov\|\cdot\|\hat x-\ox\|}{\|\hat x-\ox\|+\|\hat v-\ov\|},
$$
which implies by passing to the upper limit as $(\hat x,\hat v)\st{\tiny\gph\partial\ph}{\to}(\ox,\ov)$ that $(w,u)\in\Hat{N}((\ox,\ov);\small\gph\partial\ph)$ and thus verifies the inclusion ``$\supset$" in (\ref{pncg}), i.e., the equality therein. $\h$\vspace*{0.05in}

Now we are ready to establish the main result of this section, which is a consequence of the obtained Lemmas~\ref{polnc} and \ref{prenormal}. This result gives us the precise calculations of the prenormal cone to the subdifferential graph of an arbitrary CPWL function and hence its precoderivative entirely in terms of the {\em initial CPWL parameters} in (\ref{eq00}) and (\ref{dom}). To formulate the result and other statements below, we need the following notation. Given any index subsets $P,Q$ with $P_1\subset Q_1\subset T_1$ and $P_2\subset Q_2\subset T_2$, denote
\begin{eqnarray}\label{eq080}
\begin{array}{lll}
{\cal F}_{\tiny\{P_1,Q_1\},\{P_2,Q_2\}}:&=\span\Big\{a_i-a_j\Big|\;i,j\in P_1\Big\}+\cone\Big\{a_i-a_j\Big|\;(i,j)\in(Q_1\setminus P_1)\times P_1\Big\}\\
&+\span\Big\{d_i\Big|\;i\in P_2\Big\}+\cone\Big\{d_i\Big|\;i\in Q_2\setminus P_2\Big\},
\end{array}
\end{eqnarray}
\begin{eqnarray}\label{eq081}
\begin{array}{ll}
{\cal G}_{\tiny\{P_1,Q_1\},\{P_2,Q_2\}}:=\Big\{u\in\R^n\Big|&\la a_i-a_j,u\ra=0\;\mbox{ if }\;i,j\in P_1,\\
&\la a_i-a_j,u\ra\le 0\;\mbox{ if }\;(i,j)\in(Q_1\setminus P_1)\times P_1,\\
&\la d_i,u\ra=0\;\mbox{ if }\;i\in P_2,\;\mbox{ and }\;\la d_i,u\ra\le 0\;\mbox{ if }\;i\in Q_2\setminus P_2\;\Big\}.
\end{array}
\end{eqnarray}
Observing that the sets ${\cal F}$ and ${\cal G}$ are cones, we conclude from the classical Farkas lemma that
\begin{equation}\label{eq082}
{\cal G}_{\tiny\{P_1,Q_1\},\{P_2,Q_2\}}^*={\cal F}_{\tiny\{P_1,Q_1\},\{P_2,Q_2\}}\;\mbox{ for any }\;P_1\subset Q_1\subset T_1\;\mbox{ and }\;P_2\subset Q_2\subset T_2.
\end{equation}

Note that the following theorem presents the calculation of the aforementioned second-order constructions via the {\em fixed} index sets $P_i$ and $Q_i$, $i=1,2$, in (\ref{eq080}) and (\ref{eq081}).

\begin{Theorem}{\bf(prenormals to subdifferential graphs and precoderivatives of CPWL functions via their initial data).}\label{intial} In the setting of Lemma~{\rm\ref{polnc}} we have the following formulas with the notation {\rm(\ref{eq080})} and {\rm(\ref{eq081})}:

{\bf(i)} The prenormal cone to the subdifferential graph is calculated by
\begin{equation}\label{initial2}
 \Hat{N}((\ox,\ov);\gph\partial\ph)={\cal F}_{\tiny\{J_1,K\},\{J_2,I\}}\times{\cal G}_{\tiny\{J_1,K\},\{J_2,I\}}.
\end{equation}

{\bf(ii)} The domain and values of the precoderivative of $\ph$ at $(\ox,\ov)\in\gph\sub\th$ are calculated by
\begin{equation}\label{initial3}
(\Hat D^*\sub\ph)(\ox,\ov)(u)={\cal F}_{\tiny\{J_1,K\},\{J_2,I\}}\;\mbox{ for any }\;u\in\dom(\Hat D^*\sub\ph)(\ox,\ov)=-{\cal G}_{\tiny\{J_1,K\},\{J_2,I\}}.
\end{equation}
\end{Theorem}
{\bf Proof.} It follows from (\ref{eq074}) that $T(\ov;\partial\ph(\ox))={\cal F}_{\tiny\{J_1,K\},\{J_2,I\}}$ and from (\ref{eq060}) that $N(\ov;\partial\ph(\ox))={\cal G}_{\tiny \{J_1,K\},\{J_2,I\}}$. Thus (\ref{initial2}) is a direct consequence of Lemma~\ref{prenormal}. The precoderivative calculations in (\ref{initial3}) follow from (\ref{initial2}) due to definition (\ref{2.7}). $\h$\vspace*{0.05in}

A natural question arises on the dependence of the prenormal cone representation (\ref{initial2}), and hence the precoderivative one (\ref{initial3}), on the choice of vectors $\ov_1,\ov_2$ and multipliers $\bar\lm_i,\bar\mu_i$ in (\ref{eq04}) and (\ref{eq06}), which are not generally be unique and may potentially influence the cones ${\cal F}_{\tiny\{J_1,K\},\{J_2,I\}}$ and ${\cal G}_{\tiny\{J_1,K\},\{J_2,I\}}$ through the index sets $J_1=J_+(\ox,\ov_1)$ and $J_2=J_+(\ox,\ov_2)$ from (\ref{eq05}). The next proposition shows that it is not the case.

\begin{Proposition}{\bf (invariance of the prenormal cone and precoderivative representations for CPWL functions).}\label{invar} In the framework of Theorem~{\rm\ref{intial}} we have that the representations in {\rm(\ref{initial2})} and {\rm(\ref{initial3})} are invariant with respect to any choice of vectors $\ov_1,\ov_2$ in {\rm(\ref{eq04})}, multipliers $\bar\lm_i,\bar\mu_i$ in {\rm(\ref{eq06})}, and the index sets $J_1=J_+(\ox,\ov_1),J_2=J_+(\ox,\ov_2)$ in {\rm(\ref{eq05})}.
\end{Proposition}
{\bf Proof.}
We only need to verify the invariance of the representations in Theorem~\ref{intial} with respect to the choice of $v_1\in\co\{a_i|\;i\in K(\ox)\}$ and $v_2\in N(\ox;\dom\th)$ satisfying $\ov=v_1+v_2$. Suppose that $\ov=v_1+v_2$ with $v_1\in\co\{a_i|\;i\in K(\ox)\}$ and $v_2\in N(\ox;\dom\th)$ and then get
\begin{equation}\label{eq06002}
\begin{array}{ll}
\disp{v_1=\sum_{i\in K(\ox)}\lm'_ia_i},&\mbox{with}\;\;\disp{\sum_{i\in K(\ox)}{\lm'_i}=1}\;\mbox{ and }\;\lm'_i\ge 0,\\
\disp{v_2=\sum_{i\in I(\ox)}\mu'_id_i}&\mbox{with }\;\mu'_i\ge 0.
\end{array}
\end{equation}
Define the index set of positive multipliers in (\ref{eq06002}) by
\begin{equation}\label{eq05252}
J'_1:=\Big\{i\in K(\ox)\Big|\;\lm'_i>0\Big\}\;\mbox{ and }\;J'_2:=\Big\{i\in I(\ox)\Big|\;\mu'_i>0\Big\}.
\end{equation}
and then clarify the following relationships:
\begin{equation}\label{eq7992}
{\cal G}_{\tiny \{J_1,K\},\{J_2,I\}}={\cal G}_{\tiny\{J'_1,K\},\{J'_2,I\}}\;\mbox{ and }\;{\cal F}_{\tiny\{J_1,K\},\{J_2,I\}}={\cal F}_{\tiny\{J'_1,K\},\{J'_2,I\}}.
\end{equation}
To proceed, pick $u\in{\cal G}_{\tiny\{J_1,K\},\{J_2,I\}}$ and get for $\ov_2$ and $v$ from above that
\begin{equation}\label{eq140}
\la\ov_2,u\ra=0\;\mbox{ and }\;\la v_2,u\ra\le 0.
\end{equation}
Taking now $s\in J_1$ implies via (\ref{eq140}) that
$$
\begin{array}{ll}
\la a_s,u\ra&=\la \ov_1,u\ra=\la v_1+v_2-\ov_2,u\ra \leq \la v_1,u\ra\\
&=\Big\la\disp{\sum_{i\in K(\ox)}\lm'_ia_i,u\Big\ra=\sum_{i\in J'_1\setminus J_1}\lm'_i\la a_i,u\ra+\sum_{i\in J'_1\cap J_1}\lm'_i\la a_i,u\ra}\\
&=\disp{\sum_{i\in J'_1\setminus J_1}\lm'_i\la a_i,u\ra+\Big(\sum_{i\in J'_1\cap J_1}\lm'_i\Big)\la a_s,u\ra}.
\end{array}
$$
In this way we arrive at the estimate
\begin{equation}\label{eq8002}
\sum_{i\in J'_1\setminus J_1}\lm'_i\la a_i,u\ra\geq  \Big (1-\sum_{i\in J'_1\cap J_1} \lm'_i\Big )\la a_s,u\ra.
\end{equation}
It follows from  the equality  $1-\sum_{i\in J'_1\cap J_1}\lm'_i=\sum_{i\in J'_1\setminus J_1}\lm'_i$ along with $u\in {\cal G}_{\tiny \{J_1,K\},\{J_2,I\}}$ that
\begin{equation}\label{eq141}
{(1-\sum_{i\in J'_1\cap J_1}\lm'_i)\la a_s,u\ra=\Big (\sum_{i\in J'_1\setminus J_1}\lm'_i \Big )\la a_s,u\ra\geq \sum_{i\in J'_1\setminus J_1}\lm'_i\la a_i,u\ra}
\end{equation}
Employing next (\ref{eq8002}) together with (\ref{eq141}), we deduce that
\begin{equation}\label{eq142}
\sum_{i\in J'_1\setminus J_1}\lm'_i\la a_i,u\ra=\Big (1-\sum_{i\in J'_1\cap J_1}\lm'_i\Big )\la a_s,u\ra,
\end{equation}
which tells us that $\la a_i,u\ra=\la a_s,u\ra$ for all $i\in J'_1\setminus J_1$. Indeed, supposing $\la a_i,u\ra<\la a_s,u\ra$ for some  $i\in J'_1\setminus J_1$, leads us to a contradiction with (\ref{eq142}). This allows us to obtain
\begin{equation}\label{eq8012}
\la a_i- a_j,u\ra=0\;\mbox{ whenever }\;i,j\in J'_1 .
\end{equation}
Taking now $t\in K\setminus J'_1$ and $s\in J_1$, we claim that $\la a_t-a_s,u\ra\le 0$. Indeed, for $t\not\in J_1$ the latter inequality comes from $u\in{\cal G}_{\tiny \{J_1,K\},\{J_2,I\}}$. The opposite case of $t\in J_1$ yields $\la a_t-a_s,u\ra=0$, which therefore justifies the claim. Since $\la a_i,u\ra=\la a_s,u\ra$ for all $i\in J'_1\setminus J_1$, we get
\begin{equation}\label{eq8022}
\la a_t- a_i,u\ra\le 0\;\mbox{ whenever }\;(t,i)\in(K\setminus J'_1)\times J'_1 .
\end{equation}
It follows from the above arguments that $\la v_1-\ov_1,u\ra= 0$, and hence
$$
\la v_2,u\ra=\la\ov_2+\ov_1-v_1,u\ra=\la\ov_2,u\ra=0.
$$
The obtained relationships allow us to deduce that
$$
\la d_i,u\ra=0\;\;{\mbox{for all}}\;\;i\in J'_2\;\;\mbox{and}\;\;\la d_i,u\ra\le0\;\;\mbox{for all}\;\;i\in I\setminus J'_2,
$$
which implies together with (\ref{eq8012}) and (\ref{eq8022}) that $u\in{\cal G}_{\tiny\{J'_1,K\},\{J'_2,I\}}$, and so ${\cal G}_{\tiny\{J_1,K\},\{J_2,I\}}\subset{\cal G}_{\tiny\{J'_1,K\},\{J'_2,I\}}$. The opposite inclusion can be verified by the same arguments. Finally, the second equality in (\ref{eq7992}) follows from the first one by using the polarity in (\ref{eq082}). $\h$

\section{Second-Order Subdifferential of Piecewise Linear Functions}\label{2nd-sec}\sce

In this section we present our major calculations concerning the second-order subdifferential (\ref{secseb}) for the general class of CPWL functions. The final
formulas obtained here give us precise expressions of the domain and values of the second-order subdifferential mapping entirely via the CPWL data in (\ref{eq00}), (\ref{dom}). Several consequences of the main results are also derived below.

To begin with, we introduce some notation in addition to those formulated in Section~\ref{p-lin}. Given a pair $(\ox,\ov)\in\gph\ph$ for a CPWL function $\ph\colon\R^n\to\oR$, define the index family
\begin{equation}\label{eq090}
D(\ox,\ov):=\Big\{(P_1,P_2)\subset K(\ox)\times I(\ox)\Big|\;\ov\in\co\{a_i|\; i\in P_1\}+\cone\{d_i|\;i\in P_2\}\Big\},
\end{equation}
for any $Q_1\subset K(\ox)$ and $Q_2\subset I(\ox)$ consider the set
\begin{equation}\label{eq091}
H_{\tiny\{Q_1,Q_2\}}:=\Big\{x\in\dom\ph\Big|\;K(x)=Q_1,\;I(x)=Q_2\Big\},
\end{equation}
and then introduce the collection of index quadruples
\begin{eqnarray}\label{eq092}
\begin{array}{ll}
{\cal A}:=\Big\{(P_1,Q_1,P_2,Q_2)\Big|&P_1\subset Q_1\subset K(\ox),\;P_2\subset Q_2\subset I(\ox),\\
&(P_1,P_2)\in D(\ox,\ov),\;H_{\tiny\{Q_1,Q_2\}}\ne\emp\Big\}.
\end{array}
\end{eqnarray}
Recalling the constructions of ${\cal F}_{\tiny\{P_1,Q_1\},\{P_2,Q_2\}}$ and ${\cal G}_{\tiny\{P_1,Q_1\},\{P_2,Q_2\}}$ in (\ref{eq080}) and (\ref{eq081}), we first derive the following representation of the second-order subdifferential {\em values} for CPWL functions involving the index quadruples from (\ref{eq092}). It essentially extends the result and technique from \cite[Theorem~4.1]{hmn} developed for indicator functions of convex polyhedra.

\begin{Theorem}{\bf(second-order subdifferential values via the union of index subsets).}\label{lcplw} Let $\ph\in CPWL$ with $(\ox,\ov)\in\gph\partial\ph$, and let the set ${\cal A}$ be defined in {\rm(\ref{eq092})}. Then the limiting normal cone to $\gph\sub\ph$ at $(\ox,\ov)$ is represented as follows:
\begin{equation}\label{eq093}
N((\ox,\ov);\gph\partial\ph)=\bigcup_{(P_1,Q_1,P_2,Q_2)\in{\cal A}}{\cal F}_{\tiny\{P_1,Q_1\},\{P_2,Q_2\}}\times{\cal G}_{\tiny\{P_1,Q_1\},\{P_2,Q_2\}}.
\end{equation}
Hence the values of the second-order subdifferential $\partial^2\ph(\ox,\ov)$ at any $u\in\R^n$ are given by
\begin{eqnarray*}
\partial^2\ph(\ox,\ov)(u)=\Big\{w\in\R^n\Big|\;(w,-u)\in{\cal F}_{\tiny\{P_1,Q_1\},\{P_2,Q_2\}}\times{\cal G}_{\tiny\{P_1,Q_1\},\{P_2,Q_2\}},\;(P_1,Q_1,P_2,Q_2)\in{\cal A}\Big\}.
\end{eqnarray*}
\end{Theorem}
{\bf Proof.} It suffices to verify the normal cone representation (\ref{eq093}), which yields the second-order subdifferential formula directly by definitions (\ref{2.8}) and (\ref{secseb}). Starting with the proof of the inclusion ``$\subset$" in (\ref{eq093}), pick any $(u,w)\in N((\ox,\ov);\gph\partial\ph)$ and find by (2.2) sequences $(u_k,w_k)\to(u,w)$ and $(x_k,v_k){\st{\scriptsize\gph\sub\ph}\lra}(\ox,\ov)$ satisfying
\begin{equation}\label{eq094}
(u_k,w_k)\in\Hat N((x_k,v_k);\gph\partial\ph)\;\mbox{ for all }\;k\in\N.
\end{equation}
It follows from (\ref{eq094}) and the index set definitions in (\ref{active}), (\ref{active2}) with $I=I(\ox)$ and $K=K(\ox)$ that
$I(x_k)\subset I$ and $K(x_k)\subset K$. By passing to a subsequence if necessary, we have
\begin{eqnarray}\label{eq095}
K(x_k)=Q_1\;\mbox{ and }\;I(x_k)=Q_2\;\mbox{ whenever }\;k\in\N
\end{eqnarray}
for some index subsets $Q_1\subset I$ and $Q_2\subset K$. It implies by (\ref{eq091}) that $x_k\in H_{\tiny\{Q_1,Q_2\}}\ne\emp$. We get by (\ref{eq094}) that $v_k\in\sub\ph(x_k)$, and so Proposition~\ref{fop}(ii) tells us that $v_k=v_{1k}+v_{2k}$ for some $v_{1k}\in\co\{a_i|\;i\in K(x_k)\}$ and $v_{2k}\in N(x_k;\dom\ph)$, $k\in\N$. This allows us to show similarly to (\ref{eq06}) that there exist multipliers $\lm_{ik}\ge 0$ and $\mu_{ik}\ge 0$ for all $k\in\N$ such that
\begin{eqnarray*}
v_{1k}=\disp\sum_{i\in K(x_k)}\lm_{ik}a_i\;\mbox{ and }\;v_{2k}=\disp\sum_{i\in I(x_k)}\mu_{ik}\;\mbox{ with }\;\disp\sum_{i\in K(x_k)}\lm_{ik}=1.
\end{eqnarray*}
Extracting another subsequence of $x_k$ if needed, find index sets $P_1\subset Q_1$ and $P_2\subset Q_2$ for which
\begin{equation}\label{eq097}
\lm_{ik}>0\;\mbox{ if }\;i\in P_1,\;\lm_{ik}=0\;\mbox{ if }\;i\in Q_1\setminus P_1,\;\mu_{ik}>0\;\mbox{ if }\;i\in P_2,\;\mu_{ik}=0\;\mbox{ if }\;i\in Q_2\setminus P_2.
\end{equation}
This yields $J_+(x_k,v_{1k})=P_1$, $J_+(x_k,v_{2k})=P_2$ and thus leads us by Theorem~\ref{intial}(i) to
\begin{equation}\label{eq098}
\Hat{N}((x_k,v_k);\gph\partial\ph)={\cal F}_{\tiny\{P_1,Q_1\},\{P_2,Q_2\}}\times{\cal G}_{\tiny\{P_1,Q_1\},\{P_2,Q_2\}},\quad k\in\N,
\end{equation}
for the index sets $P_i$ and $Q_i$, $i=1,2$, from (\ref{eq095}) and (\ref{eq097}), respectively. It follows from the above and the construction of ${\cal A}$
that
$$
(u_k,w_k)\in{\cal F}_{\tiny\{P_1,Q_1\},\{P_2,Q_2\}}\times{\cal G}_{\tiny\{P_1,Q_1\},\{P_2,Q_2\}}\;\mbox{ with }\;(P_1,Q_1,P_2,Q_2)\in{\cal A},\;k\in\N,
$$
which justifies the inclusion ``$\subset$" in (\ref{eq093}) by  passing to the limit as $k\to\infty$ due to (\ref{eq098}).

To verify the opposite inclusion ``$\supset$" in (\ref{eq093}), pick a pair
$$
(u,w)\in\bigcup_{(P_1,Q_1,P_2,Q_2)\in{\cal A}}{\cal F}_{\tiny\{P_1,Q_1\},\{P_2,Q_2\}}\times{\cal G}_{\tiny\{P_1,Q_1\},\{P_2,Q_2\}}
$$
and find an index quadruple $(P_1,Q_1,P_2,Q_2)\in{\cal A}$ so that
\begin{equation}\label{eq100}
(u,w)\in{\cal F}_{\tiny\{P_1,Q_1\},\{P_2,Q_2\}}\times{\cal G}_{\tiny \{P_1,Q_1\},\{P_2,Q_2\}}.
\end{equation}
Since $H_{\tiny\{Q_1,Q_2\}}\ne\emp$ by (\ref{eq092}), find $\hat x\in H_{\tiny\{Q_1,Q_2\}}$ and get by (\ref{eq091}) that $K(\hat x)=Q_1$ and $I(\hat x)=Q_2$. Define further the convergent sequence $x_k\to\ox$ by
\begin{equation}\label{xn}
x_k=k^{-1}\hat x+(1-k^{-1})\ox,\quad k\in\N,
\end{equation}
and observe similarly to the discussions above that
\begin{equation}\label{eq101}
I(\hat x)=I(x_k)\subset I(\ox)\;\mbox{ and }\;K(\hat x)\subset K(x_k)\subset K(\ox)\;\mbox{ for large }\;k.
\end{equation}
Without loss of generality, find a constant set $Q'_1\subset K(\ox)$ such that $Q_1\subset Q'_1$ and $K(x_k)=Q'_1$ for all $k\in\N$. Let us next justify the following result.\vspace*{0.02in}

{\bf Claim:} {\em We have the equality $Q_1=Q'_1$}.\\[1ex]
It follows from (\ref{eq101}) that $Q_1\subset Q'_1$ and thus it remains to verify the opposite inclusion. Arguing by contradiction, suppose that there is $t\in Q'_1 $ with $t\not\in Q_1$. This means by construction that $x_k\in C_t$ for all $k\in\N$ while $\hat x\not\in C_t$. Define the nonempty set
\begin{eqnarray*}
E_t:=\Big\{s\in T_1\Big|\;\la a_t,\hat x\ra-\al_t<\la a_s,\hat x\ra-\al_s\;\Big\}
\end{eqnarray*}
and show that $E_t\cap K(\ox)\ne\emp$. Indeed, assuming the contrary tells us that $s\not\in K(\ox)$ for any fixed $s\in E_t$. Select now $j\in K(\hat x)$ such that $\hat x\in C_j$ and thus get $\la a_s,\hat x\ra-\al_s<\la a_j,\hat x\ra-\al_j$. Combining it with $s\in E_t$ gives us the strict inequality
\begin{eqnarray*}
\la a_t,\hat x\ra-\al_t<\la a_j,\hat x\ra-\al_j,
\end{eqnarray*}
which contradicts $\la a_j,\hat x\ra-\al_j\le\la a_t,\hat x\ra -\al_t$ due to $j\not\in E_t$ and hence justifies $E_t\cap K(\ox)\ne\emp$.

Select now $s\in E_t\cap K(\ox)$ and deduce from $t\in K(\ox)$ and (\ref{eq03}) that
\begin{equation}\label{eq1219}
\la a_t,\ox\ra-\al_t=\la a_s,\ox\ra-\al_s.
\end{equation}
Since $t\in Q'_1$, we have $x_k\in C_t$ for all $k\in\N$ and thus arrive at the inequality $\la a_s,x_k\ra-\al_s\le\la a_t,x_k\ra -\al_t$, which implies in turn that
\begin{eqnarray*}
k^{-1}[\la a_s,\hat x\ra-\al_s]+(1-k^{-1})[\la a_s,\ox\ra-\al_s]\le k^{-1}[\la a_t,\hat x\ra-\al_t]+(1-k^{-1})[\la a_t,\ox\ra-\al_t].
\end{eqnarray*}
Combining it with (\ref{eq1219}) yields $\la a_s,\hat x\ra-\al_s\le\la a_t,\hat x\ra-\al_t$. This contradicts the above choice of $s\in E_t$ and therefore verifies the claim.\vspace*{0.03in}

To continue the proof of the theorem, we get from this claim due to Proposition~\ref{fop}(ii) and formula (\ref{norc}) for $x_k$ together with (\ref{eq101}) that
\begin{equation}\label{eq103}
\sub\ph(x_k)=\co\Big\{a_i\Big|\;i\in Q_1\Big\}+\cone\Big\{d_i\Big|\;i\in Q_2\Big\}.
\end{equation}
Furthermore, it follows from (\ref{eq100}) and (\ref{eq092}) that $(P_1,P_2)\in D(\ox,\ov)$, which gives us by (\ref{eq090}) that
\begin{eqnarray}\label{Pi}
\ov=\disp{\sum_{i\in P_1}\lm_{i}a_i+\sum_{i\in P_2}\mu_{i}d_i}\;\mbox{ with some }\;\lm_i,\mu_i\ge 0\;\mbox{ and }\;\disp{\sum_{i\in P_1}\lm_{i}=1}.
\end{eqnarray}
Having the multipliers $\lm_i$ from (\ref{Pi}), we represent $P_1=P'_1\cup P''_1$ via the index subsets
$$
P'_{1}:=\Big\{i\in P_1\Big|\;\lm_i>0\Big\}\;\mbox{ and }\;P''_{1}:=\Big\{i\in P_1\Big|\;\lm_i=0\Big\}.
$$
Fix any $t\in P'_1$ and take large $k\in\N$  so that $(2k)^{-1}<\lm_t$; then for any $i\in P''_1$ select
$0<\eta_i<1$ with $\sum_{i\in P''_1}{\eta_i}=\frac{1}{2}$. Define the sequences of multipliers
$$
\lm'_{ik}:=\left\{\begin{array}{ll}
\lm_t-(2k)^{-1}&\mbox{if }\;i=t\\
k^{-1}\eta_i&\mbox{if }\;i\in P''_1\\
\lm_i&\mbox{if }\;i\in P'_1\setminus\{t\}
\end{array}\right.\quad\mbox{and}\quad\mu'_{ik}:=\mu_i+k^{-1}
$$
while observing that $\lm'_{ik}$ and $\mu'_{ik}$ are positive for $i\in P_1$ and $i\in P_2$, respectively, with
\begin{equation}\label{eq104}
\disp{\sum_{i\in P_1}\lm'_{ik}=1},\;\lm'_{ik}\to\lm_i\;\mbox{ and }\;\mu'_{ik}\to\mu_i\;\mbox{ as }\;k\to\infty,
\end{equation}
where $\lm_i$ and $\mu_i$ are taken from (\ref{Pi}). Defining further the sequences
\begin{equation}\label{vn}
v_k:=v_{1k}+v_{2k}\;\mbox{ with }\;v_{1k}:=\sum_{i\in P_1}\lm'_{ik}a_i\;\mbox{ and }\;v_{2k}:=\sum_{i\in P_2}\mu'_{ik}d_i,
\end{equation}
we deduce from (\ref{eq104}) that $v_k\to\ov$ as $k\to\infty$. It follows from (\ref{eq103}) that $v_k\in\sub\ph(x_k)$ for all $k$ due to $P_1\subset Q_1$ and $P_2\subset Q_2$. Furthermore, the positivity of $\lm'_{ik},\mu'_{ik}$ in (\ref{vn}) yields $J_+(x_k,v_{1k})=P_1$ and $J_+(x_k,v_{2k})=P_2$ while showing therefore  that
\begin{eqnarray}\label{reg-gph}
\Hat N((x_k,v_k);\gph\ph)={\cal F}_{\tiny\{P_1,Q_1\},\{P_2,Q_2\}}\times{\cal G}_{\tiny\{P_1,Q_1\},\{P_2,Q_2\}},\quad k\in\N.
\end{eqnarray}
It tells us by (\ref{eq100}) that $(u,w)\in\Hat N((x_k,v_k);\gph\ph)$ for all $k$ and finally verifies the inclusion $(u,w)\in N((\ox,\ov);\gph\ph)$ by passing to the limit in (\ref{reg-gph}) as $k\to\infty$. This justifies the normal cone representation (\ref{eq093}) and thus completes the proof of the theorem. $\h$\vspace*{0.05in}

Although the formulas of Theorem~\ref{lcplw} do not involve any ``foreign" objects for CPWL functions, they may not be so easy to get implemented while including the union ${\cal A}$ of all the index subset quadruples; cf.\ the calculations of Theorem~\ref{intial} for the `prenormal cone and precoderivative values involving only the reference index sets $I,K,J_1,J_2$. Our next theorem presents a precise calculation of the second-order subdifferential {\em domain}
$$
\dom\sub^2\ph(\ox,\ov):=\Big\{u\in\R^n\Big|\;\sub^2\ph(\ox,\ov)(u)\ne\emp\Big\}
$$
given in terms of the aforementioned index sets calculated at the reference points, i.e., entirely via the initial data. Besides being of its own importance, the obtained result is crucial for the subsequent explicit calculations of the second-order subdifferential values. It significantly extends the previous results in this direction derived in \cite{eh,hmn,hr} for the cases of component maximum functions and indicator functions of convex polyhedra, while even in these cases we arrive at enhanced formulations with a new device; see more discussions below.\vspace*{0.03in}

To proceed, fix $(\ox,\ov)\in\gph\sub\ph$, recall the notation $K:=K(\ox)$, $I:=I(\ox)$, $J_1:=J_+(\ox,\ov_1)$, and $J_2:=J_+(\ox,\ov_2)$ from (\ref{active2}), (\ref{active}), and (\ref{eq05}), and then introduce the {\em feature index subsets}
\begin{equation}\label{feature}
\begin{array}{ll}
\Gamma(J_1):=\Big\{i\in K\Big|\;\la a_i-a_j,u\ra=0\;\mbox{ for all }\;j\in J_1\;\mbox{ and }\;u\in{\cal G}_{\tiny\{J_1,K\},\{J_2,I\}}\;\Big\},\\
\Gamma(J_2):=\Big\{t\in I\Big|\;\la d_t,u\ra=0\;\mbox{ for all }\;u\in{\cal G}_{\tiny\{J_1,K\},\{J_2,I\}}\;\Big\}
\end{array}
\end{equation}
defined via ${\cal G}_{\tiny\{J_1,K\},\{J_2,I\}}$ from (\ref{eq081}) and dependent only on $(\ox,\ov)$ and the CPWL data, being invariant with respect to the choice of $\ov_1,\ov_2$ and multipliers in (\ref{eq04}); see Proposition~\ref{invar}.

\begin{Theorem}{\bf(second-order subdifferential domain for CPWL functions via initial data).}\label{domso} Let $\ph\in CPWL$ with $(\ox,\ov)\in\gph\sub\ph$. Then we have in the notation above that
\begin{equation}\label{domcod}
\dom\sub^2\ph(\ox,\ov)=\Big\{u\Big|\;\la a_i-a_j,u\ra=0\;\mbox{ for }\;i,j\in\Gamma(J_1)\;\mbox{ and }\;\la d_t,u\ra=0\;\mbox{ for }\;t\in\Gamma(J_2)\Big\}.
\end{equation}
\end{Theorem}
{\bf Proof}. To verify first the inclusion ``$\supset"$ in (\ref{domcod}), pick $u\in\R^n$ from the set on the right-hand side of (\ref{domcod}) and get $u\in{\cal G}_{\tiny\{P_1,Q_1\},\{P_2,Q_2\}}$ with $(P_1,Q_1,P_2,Q_2):=(J_1,\Gamma(J_1),J_2,\Gamma(J_2))$. Select $i\in K\setminus\Gamma(J_1)$, $t\in I\setminus\Gamma(J_2)$ and by (\ref{feature}) find $u_i\in{\cal G}_{\tiny\{J_1,K\},\{J_2,I\}}$, $j_i\in J_1$, and $u_t\in {\cal G}_{\tiny\{J_1,K\},\{J_2,I\}}$ such that
$\la a_i-a_{j_i},u_i\ra<0$ and $\la d_t,u_t\ra<0$. Define further
\begin{eqnarray}\label{ys}
y_s:=s\Big(\sum_{i\in K\setminus\Gamma(J_1)}u_i+\sum_{t\in I\setminus\Gamma(J_2)}u_t\Big)\;\mbox{ and }\;x_s:=\ox+y_s\;\mbox{ for some }\;s>0
\end{eqnarray}
and observe that $x_s-\ox=y_s\in{\cal G}_{\tiny\{J_1,K\},\{J_2,I\}}$. We split the proof into several claims.\vspace*{0.02in}

{\bf Claim~1:} {\em If $s>0$ is small enough, then $K(x_s)=Q_1$ with the above notation $Q_1:=\Gamma(J_1)$.}\\[1ex]
To verify this claim, we need to show by recalling the notation in (\ref{active2}) and (\ref{feature}) that
$$
\la a_i-a_{j},x_s\ra=\al_i-\al_j\;\mbox{ for }\;i,j\in Q_1\;\mbox{ and }\;\la a_i-a_{j},x_s\ra<\al_i-\al_j\;\mbox{ for }\;(i,j)\in(T_1\setminus Q_1)\times Q_1.
$$
Picking $i,j\in\Gamma(J_1)$ and $r\in J_1$ allows us to deduce from the constructions above that
$$
\la a_i-a_{j},x_s-\ox\ra=\la a_i-a_r,x_s-\ox\ra+\la a_r-a_{j},x_s-\ox\ra=\la a_i-a_r,y_s\ra +\la a_r-a_{j},y_s\ra=0,
$$
which in turn implies the equality
\begin{equation}\label{eq0453}
\la a_i-a_{j},x_s\ra=\la a_i-a_{j},\ox\ra=\al_i-\al_j\;\mbox{ whenever }\;i,j\in Q_1.
\end{equation}
Fix $i\in K\setminus\Gamma(J_1)$, choose $u_i\in{\cal G}_{\tiny\{J_1,K\},\{J_2,I\}}$ and $j_i\in J_1$, and observe that
\begin{eqnarray*}
\begin{array}{ll}
\la a_i-a_{j_i},x_s-\ox\ra&=\la a_i-a_{j_i},y_s\ra=s\la a_i-a_{j_i},u_i\ra\\
&+\Big\la a_i-a_{j_i},s\Big(\disp\sum_{r\in K\setminus\Gamma(J_1),r\ne i}u_r+\sum_{t\in I\setminus\Gamma(J_2)}u_t\Big)\Big\ra\leq s\la a_i-a_{j_i},u_i\ra<0.
\end{array}
\end{eqnarray*}
which yields $\la a_i-a_{j_i},x_s\ra<\la a_i-a_{j_i},\ox\ra=\al_i-\al_{j_i}$ due to $i,{j_i}\in K$. We now implement this together with (\ref{eq0453}) to arrive at the inequality
\begin{equation}\label{eq0455}
\la a_i-a_{p},x_s\ra<\al_i-\al_p\;\mbox{ for all }\;(i,p)\in(K\setminus Q_1)\times Q_1.
\end{equation}
Finally in the proof of this claim, we consider the case $i\in T_1\setminus K$. Select $j\in J_1$ and get
$$
\la a_i-a_{j},x_s\ra=\la a_i-a_{j},\ox\ra+\la a_i-a_{j},y_s\ra.
$$
Combining this with $\la a_i-a_{j},\ox\ra<\al_i-\al_j$ allows us to find $s>0$ small enough to have
\begin{equation}\label{eq0456}
\la a_i-a_{j},x_s\ra=\la a_i-a_{j},\ox\ra+\la a_i-a_{j},y_s\ra<\al_i-\al_j.
\end{equation}
Indeed, (\ref{eq0456}) is satisfied if the number $s$ is chosen in the interval
$$
0<s<\min_{(i,j)\in(T_1\setminus K)\times J_1}\Big\{\frac{\al_i-\al_j-\la a_i-a_{j},\ox\ra}{|\la a_i-a_{j},\sum_{i\in K\setminus \Gamma(J_1)}u_i+\sum_{t\in I\setminus\Gamma(J_2)}u_t\ra|}\Big\}.
$$
Therefore we deduce from (\ref{eq0453}) and (\ref{eq0456}) the estimate
\begin{eqnarray*}
\la a_i-a_{j},x_s\ra<\al_i-\al_j\;\mbox{ whenever }\;(i,j)\in(T_1\setminus K)\times Q_1,
\end{eqnarray*}
which shows together with (\ref{eq0453}) and (\ref{eq0455}) that $K(x_s)=Q_1$ for all $s$ sufficiently small.\vspace*{0.03in}

{\bf Claim 2:} {\em If $s>0$ is small enough, then $I(x_s)=Q_2$ with the above notation $Q_2:=\Gamma(J_2)$.}\\[1ex]
To verify this claim, we need to show by the notation above that
$$
\la d_t,x_s\ra=\beta_t\;\mbox{ for }\;t\in Q_2\;\mbox{ and }\;\la d_t,x_s\ra<\beta_t\;\mbox{ for }\;t\in T_2\setminus Q_2.
$$
To proceed, pick $t\in\Gamma(J_2)$ and get by construction that $\la d_t,x_s\ra=\la d_t,\ox\ra+\la d_t,y_s\ra=\beta_t$
due to (\ref{feature}) and $y_s\in{\cal G}_{\tiny\{J_1,K\},\{J_2,I\}}$ for $y_s$ from (\ref{ys}). Similarly to the proof of Claim~1, for $t\in I\setminus\Gamma(J_2)$ select $u_t\in{\cal G}_{\tiny\{J_1,K\},\{J_2,I\}}$ and arrive at the inequality
$$
\la d_t,y_s\ra=s\la d_t,u_t\ra+\Big\la d_t,s\Big(\sum_{i\in K\setminus\Gamma(J_1)}u_i+\sum_{r\in I\setminus\Gamma(J_2),r\ne t}u_r\Big)\Big\ra\le s\la d_t,u_t\ra<0.
$$
Implementing this and taking into account that $t\in I$ tell us that $\la d_t,x_s\ra<\la d_t,\ox\ra=\beta_t$.

Considering next the case of $t\in T_2\setminus I$, we get $\la d_t,\ox\ra<\beta_t$ and therefore conclude that
\begin{eqnarray*}
\la d_t,x_s\ra=\la d_t,\ox\ra+\la d_t,y_s\ra<\beta_t
\end{eqnarray*}
if $s>0$ is sufficiently. More precisely, it holds when $s$ is selected in the interval
$$
0<s<\min_{r\in T_2\setminus I}\Big\{\frac{\beta_r-\la d_r,\ox\ra}{|\la d_r,\sum_{i\in K\setminus\Gamma(J_1)}u_i+\sum_{t\in I\setminus\Gamma(J_2)}u_t\ra|}\;\Big\}.
$$
It follows from the above that $I(x_s)=Q_2$ for such $s>0$, which justifies the claim.\vspace*{0.03in}

Combining the results of Claims~1 and 2 shows that $H_{\{Q_1,Q_2\}}\ne\emp$ in (\ref{eq091}), which yields by Theorem~\ref{lcplw} that ${\cal F}_{\tiny \{P_1,Q_1\},\{P_2,Q_2\}}\subset\sub^2\ph(\ox,\ov)(u)$ for the selected quadruple $(P_1,Q_1,P_2,Q_2)$ from (\ref{eq092}). It gives us
$u\in\dom\sub^2\ph(\ox,\ov)$ thus justifying the inclusion ``$\supset$" in (\ref{domcod}).\vspace*{0.03in}

Let us continue the proof of the theorem with verifying the opposite inclusion ``$\subset$" in (\ref{domcod}). Picking $u\in\dom\sub^2\ph(\ox,\ov)$ and $w\in\sub^2\ph(\ox,\ov)(u)$ tells us by definition that $(w,-u)\in N((\ox,\ov),\gph\sub\ph)$. By Theorem~\ref{lcplw} we find a quadruple $(P_1,Q_1,P_2,Q_2)\in{\cal A}$ for which
\begin{equation}\label{2nd}
(w,-u)\in {\cal F}_{\tiny\{P_1,Q_1\},\{P_2,Q_2\}}\times{\cal G}_{\tiny\{P_1,Q_1\},\{P_2,Q_2\}}.
\end{equation}

{\bf Claim~3:} {\em We have $J_1\subset Q_1=K(x)$ and $J_2\subset Q_2=I(x)$ for some $x\in\dom\ph$ close to $\ox$. Furthermore, it holds $x-\ox\in{\cal G}_{\tiny \{J_1,K\},\{J_2,I\}}$}.\\[1ex]
To justify this claim, take $\hat x\in H_{\{Q_1,Q_2\}}$ and then construct $x_k$ by (\ref{xn}) and $v_k$ by (\ref{vn}). It is showed in the proof of Theorem~\ref{lcplw} that $(x_k,v_k){\st{\scriptsize\gph\sub\ph}\lra}(\ox,\ov)$, $K(x_k)=K(\hat x)=Q_1$, and $I(x_k)=I(\hat x)=Q_2$. Further, fix $k_0\in\N$ sufficiently large and define $x:=x_{k_0}$. Thus we get  $K(x)=Q_1$, and $I(x)=Q_2$. Appealing now to Theorem~\ref{gphpar}, we get $J_1\subset K(x)=Q_1$. It comes from $(P_1,P_2)\in D(\ox,\ov)$ that there are numbers $\lm_i\ge 0$ with $\sum_{i\in P_1}\lm_i=1$ and $\mu_i\ge 0$ such that $\ov=\sum_{i\in P_1}\lm_i a_i+\sum_{i\in P_2}\mu_id_i$. On the other hand, we know  from (\ref{eq06}) that $\ov=\sum_{i\in J_1}\bar\lm_i a_i+\sum_{i\in J_2} \bar\mu_id_i$. Remembering that $P_2\subset Q_2=I(x)$ and $I(x)\subset I$ allows us to deduce that $\la d_i,x-\ox\ra=0$ for all $i\in P_2$, which leads us to the equalities
\begin{eqnarray*}
\begin{array}{ll}
\disp{\sum_{i\in J_1}\bar\lm_i\la a_i,x-\ox\ra+\sum_{i\in J_2}\bar\mu_i\la d_i,x-\ox\ra }&=\disp{\la \ov,x-\ox\ra=\sum_{i\in P_1}\lm_i\la a_i,x-\ox\ra+\sum_{i\in P_2}\mu_i\la d_i,x-\ox\ra}\\
&=\disp{\sum_{i\in P_1}\lm_i\la a_i,x-\ox\ra+0=\sum_{i\in P_1}\lm_i\la a_i,x-\ox\ra}.
\end{array}
\end{eqnarray*}
It follows from the inclusions $J_1\subset K(x)$ and $P_1\subset Q_1=K(x)$ that
\begin{eqnarray*}
\la a_i,x\ra-\al_i=\la a_{j},x\ra-\al_j\;\;\mbox{for any}\;\;i\in J_1,\;j\in P_1.
\end{eqnarray*}
Using $P_1\subset K(\ox)$ and $J_1\subset K(\ox)$, we deduce that
\begin{eqnarray*}
\la a_i,\ox\ra-\al_i=\la a_{j},\ox\ra-\al_j\;\;\mbox{for any}\;\;i\in J_1,\;j\in P_1.
\end{eqnarray*}
Combining the above equalities brings us to $\sum_{i\in J_2}\bar\mu_i\la d_i,x-\ox\ra=0$. Since  $\la d_i,x\ra=\la d_i,\ox\ra=\beta_i$ by $\bar\mu_i\ge 0$, $i\in J$, we arrive at $J_2\subset Q_2=I(x)$. Finally, by the inclusions $J_1\subset Q_1=K(x)$ and $J_2\subset Q_2=I(x)$ we obtain $x-\ox\in{\cal G}_{\tiny\{J_1,K\},\{J_2,I\}}$. This completes the proof of this claim.\vspace*{0.03in}

{\bf Claim~4:} {\em We have the inclusions $\Gamma(J_1)\subset Q_1$ and $\Gamma(J_2)\subset Q_2$.}\\[1ex]
To verify them, suppose on the contrary that there is $i\in\Gamma(J_1)\setminus Q_1$. As mentioned in the proof of Claim~3, we have $x-\ox\in{\cal G}_{\tiny\{J_1,K\},\{J_2,I\}}$ for $x$ defined therein, and so $\la a_i-a_{j},x-\ox\ra=0$ for $j\in J_1$. On the other hand, the inclusion $i\not\in Q_1$ reduces by (\ref{pwlr1}) to $\la a_i-a_{j},x\ra<\al_i-\al_j$ for $j\in Q_1$. Since $J_1\subset Q_1$ and $\la a_i-a_{j},\ox\ra=\al_i-\al_j$, we arrive at $\la a_i-a_{j},x-\ox \ra<0$, a contradiction that justifies the first inclusion in the claim. The second one is proved similarly.\vspace*{0.07in}

{\bf Claim~5:} {\em If $u\in\R^n$ satisfies {\rm(\ref{2nd})}, then we have $-u\in{\cal G}_{\tiny\{J_1,\Gamma(J_1)\},\{J_2,\Gamma(J_2)\}}$.}\\[1ex]
Indeed, take $s\in J_1$, $r\in P_1$ and deduce from $J_1\subset Q_1$ and $-u\in{\cal G}_{\tiny\{P_1,Q_1\},\{P_2,Q_2\}}$ that $\la a_s-a_r,u\ra\ge 0$. Suppose now that $\la a_s-a_r,u\ra>0$. As showed in Claim~3, $\ov=\sum_{i\in P_1}\lm_i a_i+\sum_{i\in P_2}\mu_id_i$, which together with $-u\in {\cal G}_{\tiny\{P_1,Q_1\},\{P_2,Q_2\}}$ tells us that
\begin{equation}\label{eq0470}
\la\ov,u\ra=\sum_{i\in P_1}\lm_i\la a_i,u\ra+\sum_{i\in P_2}\mu_i\la d_i,u\ra=\sum_{i\in P_1}\lm_i\la a_r,u\ra=\la a_r,u\ra.
\end{equation}
On the other hand, we know that $\ov=\sum_{i\in J_1}\bar\lm_i a_i+\sum_{i\in J_2}\bar\mu_id_i$, and thus
\begin{equation}\label{eq0471}
\la\ov,u\ra=\sum_{i\in J_1}\bar\lm_i\la a_i,u\ra+\sum_{i\in J_2}\bar\mu_i\la d_i,u\ra\ge\sum_{i\in J_1}\bar\lm_i\la a_i,u\ra>\sum_{i\in J_1}\bar\lm_i\la a_r,u\ra=\la a_r,u\ra,
\end{equation}
It contradicts (\ref{eq0470}) and ensures in this way that
\begin{equation}\label{eq04700}
\la a_s-a_r,u\ra=0\;\mbox{ for all }\;s\in J_1\;\mbox{ and }\;r\in P_1.
\end{equation}
Thus for any $i,j\in J_1$ and $r\in P_1$ we obtain the equalities
\begin{equation}\label{eq0472}
\la a_i-a_{j},u\ra=\la a_i-a_r,u\ra+\la a_r-a_{j},u\ra=0.
\end{equation}
Take now $i\in\Gamma(J_1)\setminus J_1$ and $j\in J_1$. Picking $r\in P_1$ leads us to
\begin{equation}\label{eq0473}
\la a_i-a_{j},u\ra=\la a_i-a_r,u\ra+\la a_r-a_{j},u\ra=\la a_i-a_r,u\ra\ge 0
\end{equation}
by $\Gamma(J_1)\subset Q_1$ and $-u\in{\cal G}_{\tiny\{P_1,Q_1\},\{P_2,Q_2\}}$. It follows from (\ref{eq0470})--(\ref{eq04700}) that
$$
\sum_{i\in J_2}\bar\mu_i\la d_i,u\ra=\sum_{i\in J_1}\bar\lm_i\la a_r-a_i,u\ra=0,
$$
which gives us by $J_2\subset Q_2$ and $\bar\mu_i>0$ for $i\in J_2$ that $\la d_i,u\ra =0$. Also it follows from $\Gamma(J_2)\subset Q_2$  that
$\la d_i,u\ra\ge 0$ for any $i\in\Gamma(J_2)\setminus J_2$. Taking this into account together with (\ref{eq0472}) and (\ref{eq0473}), we complete the proof of the claim.\vspace*{0.03in}

{\bf Claim~6:} {\em We have the representation ${\cal G}_{\tiny\{J_1,\Gamma(J_1)\},\{J_2,\Gamma(J_2)\}}={\cal G}_{\tiny\{J_1,K\},\{J_2,I\}}- {\cal G}_{\tiny\{J_1,K\},\{J_2,I\}}$.}\\[1ex]
To verify the inclusion ``$\supset$" in this claim, pick $u\in{\cal G}_{\tiny\{J_1,K\},\{J_2,I\}}-{\cal G}_{\tiny\{J_1,K\},\{J_2,I\}}$ and get $u=u_1-u_2$ for some $u_1,u_2\in{\cal G}_{\tiny\{J_1,K\},\{J_2,I\}}$. Choosing $i\in\Gamma(J_1)$ and $j\in J_1$ gives us
$$
\la a_i-a_{j},u\ra=\la a_i-a_{j},u_1\ra-\la a_i-a_{j},u_2\ra=0
$$
due to (\ref{feature}). If $i\in\Gamma(J_2)$, then by using (\ref{feature}) again we arrive at
$$
\la d_t,u\ra=\la d_t,u_1\ra-\la d_t,u_2\ra=0,
$$
which implies in turn that $u\in{\cal G}_{\tiny\{J_1,\Gamma(J_1)\},\{J_2,\Gamma(J_2)\}}$ and thus justifies the inclusion ``$\supset$" above.

To verify next the opposite inclusion therein, take any $u\in{\cal G}_{\tiny\{J_1,\Gamma(J_1)\},\{J_2,\Gamma(J_2)\}}$ and observe that we are done if $K=\Gamma(J_1)$ and $I=\Gamma(J_2)$. Suppose now that either $K\setminus\Gamma(J_1)\ne\emp$ or $I\setminus\Gamma(J_2)\ne\emp$ and consider for definiteness that both sets $K\setminus\Gamma(J_1)$ and $I\setminus\Gamma(J_2)$ are nonempty while noting that the two other cases can be treated similarly. Pick $i\in K\setminus \Gamma(J_1)$ and $t\in I\setminus\Gamma(J_2) $ and find by (\ref{feature}) elements $u_i\in{\cal G}_{\tiny\{J_1,K\},\{J_2,I\}}$, $j_i\in J_1$, and  $u_t\in {\cal G}_{\tiny\{J_1,K\},\{J_2,I\}}$ such that $\la a_i-a_{j_i},u_i\ra<0$ and $\la d_t,u_t\ra<0$. Define further
$$
y_s:=s\Big(\sum_{i\in K\setminus\Gamma(J_1)}u_i+\sum_{t\in I\setminus \Gamma(J_2)}u_t\Big)\;\mbox{ and }\;x_s:=u+y_s\;\mbox{ for some }s>0.
$$
Since ${\cal G}_{\tiny\{J_1,K\},\{J_2,I\}}$ is a cone, it follows that $y_s\in{\cal G}_{\tiny\{J_1,K\},\{J_2,I\}}$. We now assert that the inclusion
$x_s\in{\cal G}_{\tiny\{J_1,K\},\{J_2,I\}}$ holds when $s$ (depending only on $u$) is sufficiently large. If it is true, then $u=x_s-y_s\in{\cal G}_{\tiny\{J_1,K\},\{J_2,I\}}-{\cal G}_{\tiny\{J_1,K\},\{J_2,I\}}$ , which would justify the claim.

To prove the assertion made, pick $i,j\in J_1$ and get by the choice of $u\in{\cal G}_{\tiny\{J_1,\Gamma(J_1)\},\{J_2,\Gamma(J_2)\}}$ and $y_s\in {\cal G}_{\tiny\{J_1,K\},\{J_2,I\}}$ the equalities
\begin{equation}\label{eq0475}
\la a_i-a_{j},x_s\ra=\la a_i-a_{j},u\ra+\la a_i-a_{j},y_s\ra=0.
\end{equation}
If $i\in\Gamma(J_1)\setminus J_1$ and $j\in J_1$, we are done by
\begin{equation}\label{eq0476}
\la a_i-a_{j},x_s\ra=\la a_i-a_{j},u\ra+\la a_i-a_{j},y_s\ra\le 0.
\end{equation}
It remains to examine the case of $i\in K\setminus\Gamma(J_1)$. Take $j_i\in J_1$ and $u_i\in {\cal G}_{\tiny\{J_1,K\},\{J_2,I\}}$, get $\la a_i-a_{j_i},u_i\ra<0$, and find $s$ so large that $\la a_i-a_{j_i},u\ra+s\la a_i-a_{j_i},u_i\ra<0$. Indeed, it works for
$$
s>\max_{i\in K\setminus\Gamma(J_1)}\Big\{0,-\frac{\la a_i-a_{j_i},u\ra}{\la a_i -a_{j_i},u_i\ra}\Big\}.
$$
This allows us to proceed with the relationships
\begin{eqnarray*}
\begin{array}{lll}
&\la a_i-a_{j_i},x_s\ra=\la a_i-a_{j_i},u\ra +\la a_i-a_{j_i},y_s\ra=\la a_i-a_{j_i},u\ra+s\la a_i-a_{j_i},u_i\ra\\
&+\Big\la a_i-a_{j_i},s\Big(\disp\sum_{r\in K\setminus\Gamma(J_1),r\ne i}u_r+\sum_{t\in I\setminus \Gamma(J_2)}u_t\Big)\Big\ra\le\la a_i-a_{j_i},u\ra+s\la a_i-a_{j_i},u_i\ra<0.
\end{array}
\end{eqnarray*}
Using this together with (\ref{eq0475}) implies for any $r\in J_1$ that
$$
\la a_i-a_r,x_s\ra=\la a_i-a_{j_i},x_s\ra+\la a_{j_i}-a_r,x_s\ra=\la a_i-a_{j_i},x_s\ra<0.
$$
Take now $t\in J_2$ and observe that $\la d_t,x_s\ra=\la d_t,u\ra+\la d_t,y_s\ra=0$, while for $t\in\Gamma(J_2)\setminus J_2$ we have $\la d_t,x_s\ra=\la d_t,u\ra+\la d_t,y_s\ra\le 0$. Finally, consider the case of $t\in I\setminus\Gamma(J_2)$ and pick $u_t\in{\cal G}_{\tiny\{J_1,K\},\{J_2,I\}}$. It follows from $\la d_t,u_t\ra<0$ that $\la d_t,u\ra+s\la d_t,u_t\ra<0$ when
$$
s>\max_{t\in I\setminus\Gamma(J_2)}\Big\{0,-\frac{\la d_t,u\ra }{\la d_t,u_t\ra}\Big\}.
$$
Therefore we are able to deduce from the above that
\begin{eqnarray*}
\begin{array}{lll}
\la d_t,x_s\ra&=&\la d_t,u\ra+\la d_t,y_s\ra=\la d_t,u\ra +s\Big\la d_t,u_t\ra+\la d_t,s\Big(\disp\sum_{i\in K\setminus\Gamma(J_1)}u_i+\disp\sum_{r\in I\setminus\Gamma(J_2),r\ne t}u_r\Big)\Big\ra\\
&\le &\la d_t,u\ra +s\la d_t,u_t\ra<0.
\end{array}
\end{eqnarray*}
Combining these pieces shows that $x_s\in{\cal G}_{\tiny\{J_1,K\},\{J_2,I\}}$ for large $s>0$, which justifies the claim.\vspace*{0.03in}

Having in hands the verified claims allows us conclude that the selected vector $u\in\R^n$ belongs to the set on right-hand side of (\ref{domcod}). Indeed, take
$i,j\in\Gamma(J_1)$ and deduce from Claims~5 and 6 that $-u=u_1-u_2$ for $u_1,u_2\in {\cal G}_{\tiny\{J_1,K\},\{J_2,I\}}$. Then it follows from (\ref{feature}) that $\la a_i-a_{j},u\ra=\la a_i-a_{j},u_2\ra-\la a_i-a_{j},u_1\ra=0$. If $t\in\Gamma(J_2)$, we again employ (\ref{feature}) to get $\la d_t,u\ra=\la d_t,u_2\ra-\la d_t,u_1\ra=0$, which thus completes the proof of the theorem. $\h$\vspace*{0.05in}

Now we consider specifications of Theorem~\ref{domso} in two important cases given by the summands in representation (\ref{theta}) of general CPWL functions. The first case below concerns the indicator function $\dd_\O(\cdot)$ of the polyhedral set $\O=\dom\ph$ corresponding to (\ref{theta}) without the maximum function part. The obtained representation is equivalent to the one in \cite[Proposition~4.4]{hmn} while being more simple and convenient for implementations.

\begin{Corollary}{\bf (second-order subdifferential domain for indicator functions of convex polyhedra).}\label{domind} Let $a_i=0\in\R^n$ and $\al_i\in\R$ for all $i\in T_1$ in the setting of Theorem~{\rm\ref{domso}}. Denoting $\O:=\dom\ph$, we have the second-order subdifferential formula
\begin{eqnarray*}
\dom\sub^2\dd_\O(\ox,\ov)=\Big\{u\in\R^n\Big|\;\la d_t,u\ra=0,\;t\in\Gamma(J_2)\Big\}.
\end{eqnarray*}
\end{Corollary}
{\bf Proof}. Follows immediately from (\ref{domcod}) with $a_i=0$ and $\al_i=0$ therein. $\h$\vspace*{0.05in}

Another particular case of (\ref{theta}) is the {\em maximum function}
\begin{equation}\label{theta1}
\ph_{\tiny{\rm max}}(x):=\max\Big\{\la a_1,x\ra-\alpha_1,\ldots,\la a_l,x\ra-\alpha_l\Big\},
\end{equation}
with corresponds to (\ref{theta}) with $\dom\ph=\R^n$. Its significant specification, where $a_i=e_i$ is the unit vector in $\R^n$ such that the
i$^{th}$ component of it is $1$ while the others are $0$ and where $\al_i=0\in\R$ for any $i\in T_1$, is given by the {\em component maximum function}
\begin{equation}\label{compo}
\phi(x):=\max\{x_1,\ldots,x_n\}\;\mbox{ for }\;x=(x_1,\ldots,x_n)\in\R^n.
\end{equation}

The next consequence of Theorem~\ref{domso} gives us a constructive formula for the second-order subdifferential domain of $\ph_{\tiny{\rm max}}$ and reduces to \cite[Theorem~3.1]{eh} in the case of (\ref{compo}).

\begin{Corollary}{\bf(second-order subdifferential domain for maximum functions).}\label{dommax} For $\ph_{\tiny{\rm max}}$ given in {\rm(\ref{theta1})} we have in the notation of Theorem~{\rm\ref{domso}} that
\begin{equation}\label{domcod3}
\dom\sub^2\ph_{\tiny{\rm max}}(\ox,\ov)=\Big\{u\in\R^n\Big|\;\la a_i-a_j,u\ra=0\;\;\mbox{for all}\;\;i,j\in\Gamma(J_1)\Big\}.
\end{equation}
In particular, for the component maximum function {\rm(\ref{compo})} we have
\begin{equation}\label{domcod4}
\dom\sub^2\phi(\ox,\ov)=\Big\{u=(u_1,\ldots,u_n)\in\R^n\Big|\;u_i=c\;\;\mbox{for all}\;\;i\in J_1\Big\}
\end{equation}
with some constant $c\in\R$.
\end{Corollary}
{\bf Proof}. Formula (\ref{domcod3}) follows immediately from Theorem~\ref{domso}. It yields representation (\ref{domcod4}) by observing that $\Gamma(J_1)=J_1$ if $a_i=e_i$ for any $i\in T_1$. $\h$\vspace*{0.05in}

Now we go back to the general case of CPWL functions $\ph$ from (\ref{eq00}) and derive results on the second-order subdifferential {\em values} expressed entirely in terms of the initial data of (\ref{eq00}) by using Theorem~\ref{domso}. First we derive an {\em upper estimate} of these values without any additional assumptions. To proceed, for $(\ox,\ov)\in\gph\ph$ and $u\in\dom\sub^2\ph(\ox,\ov)$ define the index sets
\begin{eqnarray}\label{charcset}
\begin{array}{ll}
I_{0,1}(u):=\Big\{i\in K(\ox)\Big|\;\la a_i-a_j,u\ra=0\;\;\mbox{for}\;\;j\in J_1\Big\},\\
I_{>,1}(u):= \Big\{i\in K(\ox)\Big|\;\la a_i-a_j,u\ra>0\;\;\mbox{for}\;\;j\in J_1\Big\},\\
I_{0,2}(u):=\Big\{t\in I(\ox)\Big|\;\la d_t,u\ra=0\Big\},\;I_{>,2}(u):=\Big\{t\in I(\ox)\Big|\;\la d_t,u\ra>0\;\;\mbox{for}\;\;j\in J_2\Big\}.
\end{array}
\end{eqnarray}

\begin{Theorem}{\bf(upper estimate for second-order subdifferential values of CPWL functions).}\label{upper estimate} Given $u\in\dom\sub^2\ph(\ox,\ov)$ in the setting of Theorem~{\rm\ref{domso}}, we have the inclusion
\begin{eqnarray}\label{upper1}
\begin{array}{ll}
\sub^2\ph(\ox,\ov)(u)&\subset\span\Big\{a_i-a_j\Big|\;i,j\in I_{0,1}(u)\Big\}+\cone\Big\{a_i-a_j\Big|\;i\in I_{>,1}(u),\;j\in I_{0,1}(u)\Big\}\\
&+\span\Big\{d_t\Big|\;t\in I_{0,2}(u)\Big\}+\cone\Big\{d_t\Big|\;t\in I_{>,2}(u)\Big\}
\end{array}
\end{eqnarray}
in terms of the initial CPWL parameters from {\rm(\ref{eq00})} and the index sets from {\rm(\ref{charcset})}.
\end{Theorem}
{\bf Proof}. By constructions (\ref{secseb}) and (\ref{2.8}) for $u\in\dom\sub^2\ph(\ox,\ov)$ we find $w\in\sub^2\ph(\ox,\ov)(u)$ such that $(w,-u)\in N((\ox,\ov),\gph\sub\ph)$. Applying Theorem~\ref{lcplw}, we find index subsets $(P_1,Q_1,P_2,Q_2)\in{\cal A}$ with $H_{\{Q_1,Q_2\}}\ne\emp$ and $(P_1,P_2)\in D(\ox,\ov)$ for which
$$
(w,-u)\in{\cal F}_{\tiny\{P_1,Q_1\},\{P_2,Q_2\}}\times{\cal G}_{\tiny\{P_1,Q_1\},\{P_2,Q_2\}}.
$$
Define $S_1:=\{i\in Q_1|\;\la a_i-a_j,u\ra=0\;\mbox{ if }\;j\in J_1\}$ and $S_2:=\{i\in Q_2|\;\la d_t,u\ra=0\}$ and deduce from equality (\ref{eq04700}) in the proof of Theorem~\ref{domso} that $\la a_i-a_j,u\ra=0$ whenever $j\in J_1$ and $i\in P_1$, which shows that $P_1\subset S_1$. We also have $P_2\subset S_2$ due to $-u\in{\cal G}_{\tiny\{P_1,Q_1\},\{P_2,Q_2\}}$. These observations together with $w\in{\cal F}_{\tiny\{P_1,Q_1\},\{P_2,Q_2\}}$ tell us that
$$
\begin{array}{ll}
w&\in\span\Big\{a_i-a_j\Big|\;i,j\in P_1\Big\}+\cone\Big\{a_i-a_j|\;(i,j)\in(Q_1\setminus P_1)\times P_1\Big\}\\
&+\span\Big\{d_i\Big|\;i\in P_2\Big\}+\cone\Big\{d_i\Big|\;i\in Q_2\setminus P_2\Big\}\\
&\subset\span\Big\{a_i-a_j\Big|\;i,j\in S_1\Big\}+\cone\Big\{a_i-a_j\Big|\;(i,j)\in(Q_1\setminus S_1)\times S_1\Big\}\\
&+\span\Big\{d_i\Big|\;i\in S_2\Big\}+\cone\Big\{d_i\Big|\;i\in Q_2\setminus S_2\Big\}.
\end{array}
$$
Then (\ref{upper1}) follows from $S_1\subset I_{0,1}(u)$, $S_2\subset I_{0,2}(u)$, $Q_1\setminus S_1\subset I_{>,1}(u)$, and $Q_2\setminus S_2\subset I_{>,2}(u)$. $\h$\vspace*{0.05in}

Next we establish effective conditions for the {\em equality} in (\ref{upper1}), i.e., for a precise formula to calculate the second-order subdifferential of a general CPWL function via its initial data. It is shown below that the following qualification condition is {\em sufficient} but {\em not necessary} for this.

\begin{Definition}{\bf (affine independence qualification condition).}\label{aicq} Given $\ph\in CPWL$ with $\ox\in\dom\ph$, we say the {\sc affine independence qualification condition} $($AIQC$)$ holds for $\ph$ at $\ox$ if for the generating vectors $a_i$ and $d_t$ indexed by $(i,j)\in K\times I$ with $I=I(\ox)$ and $K=K(\ox)$ the vectors $\{(a_i,1)\in\R^n\times\R|\;i\in K\}\cup\{(d_t,0)\in \R^n\times\R|\;t\in I\}$ are linearly independent.
\end{Definition}

It is clear that AIQC is implied by the linear independence of the vectors $\{a_i|\;i\in K(\ox)\}\cup\{d_t|\;t\in I(\ox)\}$, but not vice versa. Let us present some useful consequences of AIQC.

\begin{Proposition}{\bf (consequences of AIQC).}\label{cosq} Under the validity of AIQC for $\ph$ at $\ox$ we have $\Gamma(J_1)=J_1$ and $\Gamma(J_2)=J_2$ for the feature index sets in {\rm(\ref{feature})}.
\end{Proposition}
{\bf Proof.} Note that the assumed AIQC ensures the the existence of $(x,x_{n+1})\in\R^{n}\times\R$ such that
\begin{equation}\label{sys2}
\begin{array}{ll}
\la a_i,x\ra+x_{n+1}=0\;\mbox{ for }\;i\in J_1,\quad\la a_i,x\ra+x_{n+1}=-1\;\mbox{ for }\;i\in K(\ox)\setminus J_1,\\
\la d_i,x\ra=0\;\mbox{ for }\;i\in J_2,\quad\mbox{and }\;\la d_i,x\ra=-1\;\mbox{ for }\;i\in I(\ox)\setminus J_2.
\end{array}
\end{equation}
This implies by definition (\ref{eq081}) that $x\in{\cal G}_{\tiny\{J_1,K(\ox)\},\{J_2,I(\ox)\}}$, which yields in turn the claimed equalities $\Gamma(J_1)=J_1$ and $\Gamma(J_2)=J_2$. $\h$\vspace*{0.05in}

Now we are ready to accomplish the aforementioned {\em precise calculation} of second-order subdifferential values and, as a by-product, to establish the {\em second-order subdifferential sum rule} in a fully nonsmooth setting, which seems to be the first result of this type in the literature.

\begin{Theorem}{\bf (precise formula for second-order subdifferential values of CPWL functions under AIQC).}\label{percise} Assume in the setting of Theorem~{\rm\ref{upper estimate}} that AIQC holds for $\ph$ at $\ox$. Then we have the precise calculation formula
\begin{eqnarray}\label{upper2}
\begin{array}{ll}
\sub^2\ph(\ox,\ov)(u)&=\span\Big\{a_i-a_j\Big|i,j\in I_{0,1}(u)\Big\}+\cone\Big\{a_i-a_j\Big|i\in I_{>,1}(u),j\in I_{0,1}(u)\Big\}\\
&+\span\Big\{d_t\Big|t\in I_{0,2}(u)\Big\}+\cone\Big\{d_i\Big|\;t\in I_{>,2}(u)\Big\}\;\mbox{for}\;u\in\dom\partial^2\sub\ph(\ox,\ov).
\end{array}
\end{eqnarray}
Furthermore, the exact sum rule for the second-order subdifferential of $\ph=\ph_1+\ph_2$ holds:
\begin{eqnarray}\label{2nd-sum}
\sub^2\ph(\ox,\ov)(u)=\sub^2\ph_1(\ox,\ov_1)(u)+\sub^2\ph_2(\ox,\ov_2)(u),\quad u\in\dom\partial^2\sub\ph(\ox,\ov)
\end{eqnarray}
with $\ph_1(x):=\ph_{\tiny{\rm max}}(x)$ from {\rm(\ref{theta1})}, $\ph_2(x):=\dd(x;\dom\ph)$ while $\ov_1$ and $\ov_2$ are taken from {\rm(\ref{eq06})}.
\end{Theorem}
{\bf Proof}. The inclusion ``$\subset$" is obtained in (\ref{upper1}). To verify the opposite inclusion, pick any $w$ from the set on the right-hand side in (\ref{upper2}). Select $P_1=I_{0,1}(u)$, $Q_1=I_{0,1}(u)\cup I_{>,1}(u) $, $P_2=I_{0,2}(u)$, and $Q_2=I_{0,2}(u)\cup I_{>,2}(u)$ and get
\begin{equation}\label{eq0490}
(w,-u)\in {\cal F}_{\tiny\{P_1,Q_1\},\{P_2,Q_2\}}\times{\cal G}_{\tiny\{P_1,Q_1\},\{P_2,Q_2\}}.
\end{equation}
Define $x_t:=\ox+tx$ for $t>0$, where $x$ solves the system
\begin{equation}\label{sys3}
\begin{array}{ll}
\la a_i-a_j,x\ra=0\;\mbox{ for }\;i,j\in Q_1,\quad\la a_i-a_j,x\ra=-1\;\mbox{ for }\;(i,j)\in(K(\ox)\setminus Q_1)\times Q_1,\\
\la d_i,x\ra=0\;\mbox{ for }\;i\in Q_2,\;\mbox{ and }\;\la d_i,x\ra=-1\;\mbox{ for }\;i\in I(\ox)\setminus Q_2.
\end{array}
\end{equation}
Such a solution exists under AIQC due to Proposition~\ref{cosq}, which allows us to replace $Q_1$ by $J_1$ and $\O_2$ by $J_2$ and thus to reduce (\ref{sys3}) to system (\ref{sys2}) considered above.

We claim now that $x_t\in H_{\{Q_1,Q_2\}}$ in (\ref{eq091}). Indeed, it follows directly from the equalities $K(x_t)=Q_1$ and $I(x_t)=Q_2$ when $t>0$ is sufficiently small. To finish the proof of (\ref{upper2}), it remains to show that $(P_1,P_2)\in D(\ox,\ov)$ with $D(\ox,\ov)$ taken from (\ref{eq090}). Since $u\in \dom\sub^2\ph(\ox,\ov)$, it follows from (\ref{domcod}) that $J_1\subset I_{0,1}(u)=P_1$ and $J_2\subset I_{0,2}(u)=P_2$. On the other hand, we know from (\ref{eq06}) that $(J_1,J_2)\in D(\ox,\ov)$, which ensures therefore that $(P_1,P_2)\in D(\ox,\ov)$. The latter implies that $(P_1,Q_1,P_2,Q_2)\in{\cal A}$ with ${\cal A}$ taken from (\ref{eq092}). Employing finally Theorem~\ref{lcplw} together with (\ref{eq0490}) yields $w\in\sub^2\ph(\ox,\ov)(u)$ and thus justifies (\ref{upper2}). The second-order sum rule in (\ref{2nd-sum}) is an immediate subsequence of (\ref{upper2}). $\h$\vspace*{0.05in}

The precise calculation of the second-order subdifferential for CPWL functions in (\ref{upper2}) via their initial data opens the gate for a variety of applications to numerous issues of variational analysis and optimization, which will be considered in our subsequent research. Recall that compositions of CPWL functions with ${\cal C}^2$-mappings produce a major subclass of {\em fully amenable} functions frequently appeared in many aspects of variational analysis and optimization; see, e.g., \cite{rw} and also \cite{mr} for more recent developments. The exact chain rule for such compositions was derived in \cite[Theorem~4.3]{mr} under a certain {\em second-order qualification condition} involving the second-order subdifferential of a CPWL outer function. The precise calculation in (\ref{upper2}) makes these results more efficient for further implementations.

\section{Some Particular Cases}\label{spec}\sce

In this section we revisit the second-order subdifferential calculations for some types of the maximum functions obtained recently in \cite{eh} and \cite{e14} and then derive, being motivated by these papers, new results in this direction based on Theorem~\ref{lcplw}. Let us start with the component maximum function (\ref{compo}). The following result was obtained in \cite[Theorem~3.1]{eh} by reducing the problem to the polyhedral framework of \cite{dr96} with the subsequent usage of the critical face condition. Our approach is more direct presenting a straightforward application of Theorem~\ref{percise}.

\begin{Proposition} {\bf(calculating the second-order subdifferential of the component maximum function).}\label{subcomp} Let $\phi$ be given in {\rm(\ref{compo})} with $u=(u_1,\ldots,u_n)\in\dom\sub^2\phi(\ox,\ov)$. Then we have $u_i=\gg$ for all $i\in J_1=J_+(\ox,\ov_1)$ for some constant $\gg\in\R$ and also $I_{>,1}(u)=\{i\in K(\ox)|\;u_i>\gg\}$. Furthermore, the following representation holds:
 \begin{equation}\label{max4}
\begin{array}{ll}
\partial^2\phi(\ox,\ov)(u)=\Big\{w=(w_1,\ldots,w_n)\in\R^n\Big|&\disp\sum_{i=1}^{n}w_i=0,\;w_i\ge 0\;\mbox{ if }\;i\in I_{>,1}(u),\\
&w_i=0\;\mbox{ if }\;i\in(T_1\setminus K(\ox))\cup I_{<,1}(u)\Big\},
\end{array}
\end{equation}
where the index set $I_{<,1}(u):=\{i\in K(\ox)|\;u_i<\gg\}$.
\end{Proposition}
{\bf Proof}. Given $u\in\dom\sub^2\phi(\ox,\ov)$, it comes from (\ref{domcod4}) that $u_i=\gg$  with some $\gg\in\R$ whenever $i\in J_1$. Remembering that $a_i=e_i\in\R^n$ as $i\in T_1$ and plugging this into (\ref{charcset}) give us $I_{>,1}(u)=\{i\in K(\ox)|\;u_i>u_j\;\mbox{for }\;j\in J_1\}=\{i\in K(\ox)|\;u_i>\gg\} $. Now take $w\in\partial^2\phi(\ox,\ov)(u)$ and observe the validity of AIQC in Theorem~\ref{percise}. Thus we get
\begin{equation}\label{upper3}
\sub^2\phi(\ox,\ov)(u)=\span\Big\{e_i-e_j\Big|\;i,j\in I_{0,1}(u)\Big\}+\cone\Big\{e_i-e_j\Big|\;i\in I_{>,1}(u),\;j\in I_{0,1}(u)\;\Big\},
\end{equation}
where $I_{0,1}(u)=\{i\in K(\ox)|\;u_i=\gg\}$. This leads us to the representation
$$
w=\sum_{(i,j)\in I_{0,1}(u)\times I_{0,1}(u)}\lm_{ij}(e_i-e_j)+\sum_{(i,j)\in I_{>,1}(u)\times I_{0,1}(u)}\mu_{ij}(e_i-e_j)\;\mbox{ with }\;\lm_{ij}\in\R,\;\mu_{ij}\ge 0,
$$
which implies that $\sum_{i=1}^{n}w_i=0$ and $w_i\ge 0$ for $i\in I_{>,1}(u)$. Since $I_{0,1}(u)\subset K(\ox)$ and $I_{>,1}(u)\subset K(\ox)$, there is no vector $e_i$ with $i\in(T_1\setminus K(\ox))\cup I_{<,1}(u)$ in (\ref{upper3}). This yields $w_i=0$ for $i\in(T_1\setminus K(\ox))\cup I_{<,1}(u)$ and shows that $w$ belongs to the set on the right-hand side of (\ref{max4}). To verify the opposite inclusion, take $w$ from the latter set for some $u\in\dom\sub^2\phi(\ox,\ov)$ and get
$$
\sum_{i\in I_{>,1}(u)}w_i+\sum_{i\in I_{0,1}(u)}w_i =0\;\mbox{ with }\;w_i\ge 0\;\mbox{ for }\;i\in I_{>,1}(u).
$$
If $w=0$, then we have $w\in\partial^2\phi(\ox,\ov)(u)$ due to (\ref{upper3}). Otherwise, let $r\in I_{0,1}(u)$ and observe that
$w_r=-(\sum_{i\in I_{>,1}(u)} w_i+\sum_{i\in I_{0,1}(u),i\ne r}w_i)$, which leads us to
$$
w=\sum_{i\in I_{>,1}(u)}w_i(e_i-e_r)+\sum_{i\in I_{0,1}(u),i\ne r}w_i(e_i-e_r)
$$
and yields $w\in\sub^2\phi(\ox,\ov)(u)$ by the above representation of $w$. $\h$\vspace*{0.05in}

To proceed further, note that we Theorem~\ref{lcplw} for the maximum function (\ref{theta1}) gives us
\begin{eqnarray}\label{cor-max}
\partial^2\ph_{\tiny{\rm max}}(\ox,\ov_1)(u)=\Big\{w\in\R^n\Big|\;(w,-u)\in{\cal F}_{\tiny\{P_1,Q_1\}}\times{\cal G}_{\tiny\{P_1,Q_1\}}\;\mbox{ with }\;(P_1,Q_1)\in{\cal A}\Big\},
\end{eqnarray}
where $u\in\R^n$ and $\ov_1$ is from {\rm (\ref{eq04})}. This fact has not been observed in \cite{eh}, where the authors only consider either the case of $a_i=e_i$ for $i\in T_1$, or the case where the vectors $a_i$, $i\in K(\ox)$, are affinely independent in the sense of Definition~\ref{aicq}; see \cite[Theorem~3.1 and Theorem~4.2]{eh}. Then it was explored in \cite{e14} for some other particular cases of the maximum function under the validity of AIQC. In the rest of this section we obtain, on the basis of (\ref{cor-max}), new results in this direction that cover the aforementioned particular cases under AIQC while also encompass more general settings where AIQC does not hold and the results of \cite{eh} and \cite{e14} cannot be applied.

First we address  the so-called {\em $\infty$-norm} function
\begin{eqnarray}\label{infin}
\phi_\infty(x):=\|x\|_\infty=\max\Big\{|x_1|,\ldots,|x_n|\Big\}\;\mbox{ for }\;x=(x_1,\ldots,x_n)\in\R^n.
\end{eqnarray}
To rewrite (\ref{infin}) in form (\ref{theta1}), observe that $\phi_\infty(x)=\max\{x_1,\ldots,x_n,-x_1,\dots,-x_n\big\}$ and so
$$
\psi_\infty(x)=\max\Big\{\la a_1,x\ra,\ldots,\la a_{2n},x\ra\Big\}\;\mbox{ with }\;a_i=e_i\;\mbox{ and }\;a_{n+i}=-e_i,\quad i=1,\ldots,n.
$$
Explicit formulas for calculating $\partial^2\phi_\infty(\ox,\ov)$ entirely via the initial data were derived in \cite{e14} in the two cases: (a) $\ox=0$ with $\ov\in\int\sub\phi_\infty(\ox)$ and (b) $\ox\ne 0$. The most delicate and important case of $\ox=0$ with $\ov\in\bd\sub\phi_\infty(\ox)$ has not been resolved in \cite{e14} due to the violation of AIQC in this setting. Now we are able to proceed in this case based on Theorem~\ref{upper estimate} and formula (\ref{cor-max}).\vspace*{0.03in}

The next calculation formulas for the domain and values of $\partial^2\phi_\infty(\ox,\ov)$ make use of the classical sign function $\sgn(x)$ equal to $1$ if $x>0$, to $0$ if $x=0$, and to $-1$ if $x<0$.

\begin{Theorem}{\bf (second-order subdifferential of the $\infty$-norm function).}\label{infinity} Considering the $\infty$-norm function {\rm(\ref{infin})} with $\ox=0$ and $\ov=(\ov_1,\ldots,\ov_n)\in\bd\sub\phi_\infty(\ox)$, define
$$
J_\infty:=\Big\{i\in\{1,\ldots,n\}\Big|\;\ov_i\ne 0\Big\}\;\mbox{ and }\;J_\infty^c:=\{1,\ldots,n\}\setminus J_\infty.
$$
Then we have the following assertions:

{\bf(i)} The domain of $\sub^2\phi_\infty$ at $(\ox,\ov)$ is calculated by
$$\
\dom\sub^2\phi_\infty(\ox,\ov)=\Big\{u=(u_1,\ldots,u_n)\in\R^n\Big|\;\sgn(\ov_i)u_i=\gg\;\;\mbox{for all}\;\;i\in J_\infty\Big\},
$$
where $\gg\in\R$ is some real constant.

{\bf(ii)} We have the following upper bound for the value of $\sub^2\phi_\infty(\ox,\ov)$ at any $u\in\dom\sub^2\phi_\infty(\ox,\ov)$:
\begin{eqnarray}\label{inf-ii}
\begin{array}{ll}
\sub^2\phi_\infty(\ox,\ov)(u)&\subset\span\Big\{\nu_ie_i-\nu_je_j\Big|\;i,j\in L_1(u)\cup J_\infty\Big\}\\
&+\cone\Big\{\mu_ie_i-\nu_je_j\Big|\;(i,j)\in L_2(u)\times(L_1(u)\cup J_\infty)\Big\},
\end{array}
\end{eqnarray}
where the index sets $L_1(u)$ and $L_2(u)$ are defined by
$$
L_1(u):=\Big\{i\in J_\infty^c\Big|\;u_i=\gg\;\mbox{or}\;-u_i=\gg\Big\}\;\mbox{ and }\;L_2(u):=\Big\{i\in J_\infty^c\Big|\;u_i>\gg\;\mbox{ or }\;-u_i>\gg\Big\}
$$
with the constant $\gg\in\R$ taken from {\rm(i)} and
$$
\nu_i:=\left\{\begin{array}{ll}
1&\mbox{if}\;\;{u_i=\gg}\\
-1&\mbox{if}\;\;{-u_i=\gg}
\end{array}\right.
\quad
\mbox{for}\quad i\in L_1(u)\cup J_\infty,
$$
$$
\mu_i:=\left\{\begin{array}{ll}
1&\mbox{if}\;\;{u_i>\gg}\\
-1&\mbox{if}\;\;{-u_i>\gg}
\end{array}\right.
\quad
\mbox{for}\quad i\in L_2(u).
$$
Furthermore, equality holds in {\rm(\ref{inf-ii})} for $u=0$ and for $u=(u_1,\ldots,u_n)\in\dom\sub^2\phi_\infty(\ox,\ov)$ satisfying $\min\{u_i,-u_i\}<\gg$ whenever $i\in \{1,\ldots,n\}$.
\end{Theorem}
{\bf Proof.} It is easy to see that $K(\ox)=\{1,\ldots,2n\}$ for the active index set (\ref{active2}) when $\ox=0$. Furthermore, we have the representations
$$
\ov\in\sub\phi_\infty(\ox)=\co\Big\{a_1,\ldots,a_{2n}\Big\}\Longleftrightarrow\ox\in N(\ov;\B_1),
$$
where $\B_1:=\{y=(y_1,\ldots,y_n)|\;\|y\|_1:=\sum_{i=1}^{n}|y_i|\le 1\}$, and where $a_i=e_i\in\R^n$ and $a_{n+i}=-e_i\in\R^n$ for $i=1,\ldots,n$. Thus there are multiplies $\bar\lm_i\ge 0$ for $i\in K(\ox)$ such that
\begin{equation}\label{ki2}
\ov=\sum_{i\in K(\ox)}\bar\lm_i a_i\;\mbox{ and }\;\sum_{i\in K(\ox)}\bar\lm_i=1\;\mbox{ with }\;\bar\lm_i\bar\lm_{i+n}=0\;\mbox{for}\;i\in K(\ox)\cap \{1,\ldots,n\}.
\end{equation}
Remembering from (\ref{eq05}) that $J_1=J_+(\ox,\ov_1)=\{i\in K(\ox)|\;\bar\lm_i>0\}$ with $\ov_1=\ov$ and using $\ov\in\bd\sub\phi_\infty(\ox)$, we get that $\|\ov\|_1=\sum_{i=1}^{n}|\ov_i|=1$. Combining this with (\ref{ki2}) gives us
\begin{equation}\label{ki20}
\bar\lm_i=\left\{\begin{array}{ll}
|\ov_i|&\mbox{if }\;i\le n\\
|\ov_{i-n}|&\mbox{if }\;i>n
\end{array}\right.\quad\mbox{and}\quad a_i=\left\{\begin{array}{ll}
\sgn(\ov_i)e_i&\mbox{if }\;i\le n\\
\sgn(\ov_{i-n})e_{i-n}&\mbox{if }\;i>n
\end{array}\right.\quad
\mbox{for}\quad i\in J_1.
\end{equation}
We deduce from (\ref{ki20}) the implications
\begin{equation}\label{ki21}
\Big[i\in J_1\Longrightarrow\left\{\begin{array}{ll}
i\in J_\infty&\mbox{if}\;i\le n\\
i-n\in J_\infty&\mbox{if}\;i>n
\end{array}\right.\Big],\quad\Big[i\in J_\infty\Longrightarrow\;\left\{\begin{array}{ll}
i\in J_1&\mbox{if}\;\sgn(\ov_i)>0\\
i+n\in J_1&\mbox{if}\;\sgn(\ov_i)<0
\end{array}\right.\Big].
\end{equation}

To justify now assertion (i), we apply Theorem~\ref{dommax} along with (\ref{domcod3}) and arrive at
\begin{equation}\label{ki5}
\dom\sub^2\phi_\infty(\ox,\ov)=\Big\{u=(u_1,\ldots,u_n)\in\R^n\Big|\;\la a_i-a_j,u\ra=0\;\;\mbox{for all}\;\;i,j\in\Gamma(J_1)\Big\}.
\end{equation}
Since $\dom\phi_\infty=\R^n$, we drop for this function the indexes related to the second summand in (\ref{theta}). Letting $K:=K(\ox)$ and taking into account that $J_1\subset K$ together with (\ref{eq081}), (\ref{ki20}), and (\ref{ki21}), we prove the following claim that gives us a constructive representation of ${\cal G}_{\tiny\{J_1,K\}}$:\vspace*{0.02in}

{\bf Claim:} {\em The set ${\cal G}_{\tiny\{J_1,K\}}$ from {\rm(\ref{eq081})} in the case of {\rm(\ref{infin})} admits the representation
\begin{equation}\label{ki4}
u=(u_1,\ldots,u_n)\in{\cal G}_{\tiny\{J_1,K\}}\Longleftrightarrow\left\{\begin{array}{ll}
\sgn(\ov_i)u_i=\sgn(\ov_j)u_j&\mbox{if }\;{i,j\in J_\infty,}\\
|u_i|\le\sgn(\ov_j)u_j&\mbox{if }\;{(i,j)\in J_\infty^c\times J_\infty},\\
u_i\ge 0&\mbox{if }\;\sgn(\ov_j)>0,\\
u_i\le 0&\mbox{if }\;\sgn(\ov_j)<0.
\end{array}\right.
\end{equation}}\\[1ex]
We begin with verifying the inclusion ``$\subset$" in (\ref{ki4}). Pick $u\in{\cal G}_{\tiny\{J_1,K\}}$ and get by (\ref{eq081}) that
\begin{equation}\label{eq511}
\la a_i-a_j,u\ra=0\;\mbox{ for all }\;i,j\in J_1.
\end{equation}
Letting first $i,j\in J_\infty$ gives us $i,j\le n$, $\ov_i\ne 0$, and $\ov_j\ne 0$. Consider the following cases:

{\bf (1)} $\sgn(\ov_i)>0,\;\sgn(\ov_j)>0$. This leads us to $i,j\in J_1$, $a_i=\sgn(\ov_i)e_i$, and $a_j=\sgn(\ov_j)e_j$. We get $\sgn(\ov_i)u_i=\sgn(\ov_j)u_j$ by (\ref{eq511}), and thus $u$ belongs to the right-hand side of (\ref{ki4}).

{\bf (2)} $\sgn(\ov_i)>0,\;\sgn(\ov_j)<0$. This means that $i\in J_1$, $j+n\in J_1$, and therefore we conclude that  $a_i=\sgn(\ov_i)e_i$ and $a_j=\sgn(\ov_{j+n-n})e_{j+n-n}=\sgn(\ov_j)e_j$.
Using (\ref{eq511}) implies that $\sgn(\ov_i)u_i=\sgn(\ov_j)u_j$ and thus $u$ belongs to the right-hand side of (\ref{ki4}).

{\bf (3)} $\sgn(\ov_i)<0,\;\sgn(\ov_j)>0$. We can treat this case similarly to (2).

{\bf (4)} $\sgn(\ov_i)<0,\;\sgn(\ov_j)<0$. We can treat this case similarly to (1).

Consider next the case of $(i,j)\in J_\infty^c\times J_\infty $ in (\ref{ki4}), which gives us $j\le n$ and  $\ov_j\ne 0$. Let for definiteness $i\le n$ observing that the other case of $i>n$ can be treated similarly. Supposing that $i\in J_1$, we deduce from (\ref{ki21}) that $i\in J_\infty$, a contradiction. This justifies that $i\in K\setminus J_1$, which together with the condition $i\le n$ tells us that $a_i=e_i$. Employing now (\ref{eq081}) shows that
\begin{equation}\label{eq930}
u_i=\la a_i,u\ra\le\la a_j,u\ra=\sgn(\ov_j)e_j.
\end{equation}
On the other hand, we observe that $i+n\in K\setminus J_1$. Indeed, if $i+n\in J_1$, it follows from (\ref{ki21}) that $i=i+n -n\in J_\infty$, a contradiction. Using this observation confirms that $a_{i+n}=-e_i$, and thus we get $-u_i=\la a_{i+n},u\ra\le\la a_j u\ra=\sgn(\ov_j)e_j$ by (\ref{eq081}). Taking it into account along with (\ref{eq930}) yields $|u_i|\le\sgn(\ov_j)e_j$, and thus $u$ belongs to the right-hand side of (\ref{ki4}).

Considering further the case of $\sgn(\ov_j)>0$ in (\ref{ki4}), we need to verify that $u_i\ge 0$. To furnish this, deduce from $\sgn(\ov_j)>0$ that $j\in J_\infty$, and thus it follows from (\ref{ki21}) that $j\in J_1$. Since $j\le n$, we get from (\ref{ki20}) that $a_j=e_j$. Moreover, it results from (\ref{ki2}) that $j+n\in K\setminus J_1$, which tells us that $a_{j+n}=-e_j$. Appealing now to (\ref{eq081}) implies that
\begin{eqnarray*}
-u_j=\la a_{j+n},u\ra\le\la a_j,u\ra=u_j,
\end{eqnarray*}
and thus leads us to $u_j\ge 0$. The remaining case of $\sgn(\ov_j)<0$ in (\ref{ki4}) can be done in similarly to the the previous one, and thus we finishes the proof of the inclusion ``$\subset$" in (\ref{ki4}).\vspace*{0.01in}

Our next task is to verify the inclusion ``$\supset$" in (\ref{ki4}). Let us first show that equality (\ref{eq511}) holds if $i,j\in J_1$. We split the proof of it into the following cases:

{\bf(a)} $i,j\le n$. This yields (\ref{ki21}) by $i,j\in J_\infty$. Observe from (\ref{ki20}) that $a_i=\sgn(\ov_i)e_i$ and that $a_j=\sgn(\ov_j)e_j$. Using it along with  the relations in (\ref{ki4}) confirms the validity of (\ref{eq511}).

{\bf(b)} $i>n,\;j\le n$. Employing (\ref{ki21}) in this case leads us to $i-n\in J_\infty$ and $j\in J_\infty$. It allows us to deduce from the relations in (\ref{ki4}) and (\ref{ki20}) that
$$
\la a_i,u\ra=\la\sgn(\ov_{i-n})e_{i-n},u\ra=\sgn(\ov_{i-n})u_{i-n}=\sgn(\ov_j)u_j=\la\sgn(\ov_j)e_j,u\ra=\la a_j,u\ra,
$$
which thus verifies the validity of (\ref{eq511}) in this case.

{\bf(c)} $i\le n,\;j>n$. This case can be treated similar to case (b).

{\bf(d)} $i>n,\;j > n$. We can prove that equality (\ref{eq511}) holds similarly to case (a).

To finish the proof of the claim, it remains to justify the validity of the inequality $\la a_i-a_j,u\ra\le 0$ provided that $i\in K\setminus J_1$ and $j\in J_1$.  Here again we need to consider several different cases for $i$ and $j$. For brevity we consider only the case of $i,j\le n$ while observing that the other cases can be done similarly. It follows from the the relations in (\ref{ki21}) and (\ref{ki20}) that $j\in J_\infty$ and that $a_j=\sgn(\ov_j)e_j$. Since $i\le n$, we get $a_i=e_i$. If $i\in J_\infty^c$, then it follows from (\ref{ki4}) that
$$
\la a_i,u\ra =u_i\le\sgn(\ov_j)u_j=\la\sgn(\ov_j)e_j,u\ra=\la a_j,u\ra,
$$
which brings us to the claimed inequality. Otherwise we get $i\in J_\infty$ and observe that $\sgn(\ov_i)<0$ and hence $u_i \le 0$ by (\ref{ki4}). Considering further $i\le n$ tells us that $a_i=e_i$, and thus it follows from the relations in (\ref{ki4}) that
$$
\la a_i,u\ra=u_i\le-u_i=\sgn(\ov_i)u_i=\sgn(\ov_j)u_j=\la\sgn(\ov_j)e_j,u\ra=\la a_j,u\ra,
$$
which therefore completes the proof of the claim.\vspace{0.03in}

To continue with the proof of the theorem, pick $i\in\Gamma(J_1)$ for $j\in J_1$ and suppose that $j\le n$. It follows from (\ref{feature}) and (\ref{ki20}) that
$\la a_i,u\ra=\la a_j,u\ra=\sgn(\ov_j)u_j$ for all $u\in{\cal G}_{\tiny\{J_1,K\}}$. This together with (\ref{ki4}) yields $i\in J_1$, since otherwise there is $u$ from (\ref{ki4}) violating $\la a_i,u\ra=\sgn(\ov_j)u_j$. To see it, suppose that $i\in\Gamma(J_1)$ and $i\not\in J_1$. Suppose further that $i\le n$; the other case $i>n$ can be treated similarly. There are two possible subcases here: either (a) $i\in J_\infty$, or (b) $i\not\in J_\infty^c$. If in the first (sub)case we assume that $\sgn(\ov_i)>0$, then (\ref{ki21}) yields $i\in J_1$, a contradiction. Thus we have $\sgn(\ov_i)<0$ that leads us by (\ref{ki21}) to $i+n\in J_1$.
It follows from $u\in{\cal G}_{\tiny\{J_1,K\}}$ that $\la a_{i+n}-a_j,u\ra=0$, which tells us that $-u_i=\sgn(\ov_j)u_j$. Employing the latter together with $u_i=\la a_i,u\ra=\sgn(\ov_j)u_j$ yields $u_j=0$. This means that for any $u\in{\cal G}_{\tiny\{J_1,K\}}$ we must have $u_j=0$, which contradicts (\ref{ki4}).

Considering further (sub)case (b) above, we get $(i,j)\in J_\infty^c\times J_\infty $. Select $u\in{\cal G}_{\tiny\{J_1,K\}}$ so that $u_i<\sgn(\ov_j)u_j$, which is possible due to (\ref{ki4}). This contradicts the fact that $u_i=\la a_i,u\ra=\sgn(\ov_j)u_j$ for all $u\in{\cal G}_{\tiny\{J_1,K\}}$ and thus verifies
the inclusion $\Gamma(J_1)\subset J_1$. Similar arguments work for the case of $j>n$. Since the opposite inclusion $J_1\subset\Gamma(J_1)$ is trivial, we arrive at $\Gamma(J_1)=J_1$ and, combining it with (\ref{ki21}) and (\ref{ki5}), justify assertion (i) of the theorem.\vspace*{0.02in}

Next we verify assertion (ii). Pick $u\in\dom\sub^2\phi_\infty(\ox,\ov)$ and get from Theorem~\ref{upper estimate} that
\begin{equation}\label{ki7}
\sub^2\phi_\infty(\ox,\ov)(u)\subset\span\Big\{a_i-a_j\Big|\;i,j\in I_{0,1}(u)\Big\}+\cone\Big\{a_i-a_j\Big|\;i\in I_{>,1}(u),\;j\in I_{0,1}(u)\Big\}.
\end{equation}
Implementing (\ref{charcset}) and (\ref{ki20}) gives us the equivalences
\begin{equation}\label{ki22}
i\in I_{0,1}(u)\Longleftrightarrow i\in K(\ox),\quad\;\la a_i,u\ra=\gg\Longrightarrow\left\{\begin{array}{ll}
i\in L_1(u)\cup J_\infty&\mbox{if }\;{i\le n,}\\i-n\in L_1(u)\cup J_\infty&\mbox{if }\;i>n,\end{array}\right.
\end{equation}
\begin{equation}\label{ki228}
i\in L_1(u)\cup J_\infty \Longrightarrow\left\{\begin{array}{ll}
i\in I_{0,1}(u)&\mbox{if }\;{u_i=\gg,}\\
i+n\in I_{0,1}(u)&\mbox{if }\;{-u_i=\gg,}\\
\end{array}\right.
\end{equation}
\begin{equation}\label{ki23}
i\in I_{>,1}(u)\Longleftrightarrow i\in K(\ox),\quad\la a_i,u\ra>\gg\Longrightarrow\left\{\begin{array}{ll}
i\in L_2(u)&\mbox{if }\;i\le n,\\i-n\in L_2(u)&\mbox{if }\;i>n,
\end{array}\right.
\end{equation}
\begin{equation}\label{ki2288}
i\in L_2(u) \Longrightarrow\left\{\begin{array}{ll}
i\in I_{>,1}(u) &\mbox{if }\;{u_i>\gg,}\\
i+n\in I_{>,1}(u) &\mbox{if }\;{-u_i>\gg},\\
\end{array}\right.
\end{equation}
which justifies inclusion (\ref{inf-ii}) in (ii). To prove further the equality therein in the cases claimed in the theorem let us verify that the inclusion ``$\supset$" holds in (\ref{ki7}) when either (a) the vector $u=(u_1,\ldots,u_n)\in\dom\sub^2\phi_\infty(\ox,\ov)$ satisfies $\min\{u_i,-u_i\}<\gg$ for all $i\in \{i,\ldots,n\}$, or (b) $u=0$. To proceed, pick $w$ from the set on the right-hand side in (\ref{ki7}) and denote $P_1:=I_{0,1}(u)$ and $Q_1:=I_{0,1}(u)\cup I_{>,1}(u)$. This yields
\begin{equation}\label{ki8}
(w,-u)\in{\cal F}_{\tiny\{P_1,Q_1\}}\times{\cal G}_{\tiny\{P_1,Q_1\}}.
\end{equation}
Observe furthermore that in case (a) we have $K(x)=Q_1$ with $x=(x_1,\ldots,x_n)$ given by
$$
x_i:=\left\{\begin{array}{ll}
1&\mbox{if }\;i\in Q_1\cap\{1,\ldots,n\},\\
-1&\mbox{if }\;n+i\in Q_1\cap\{n+1,\ldots,2n\}\\
0&\mbox{otherwise}
\end{array}\right.
$$
for all $i=1,\ldots,n$. In case (b) we also have $K(x)=Q_1$ for $x=0$. This results in $H_{\{Q_1\}}\ne\emp$ for the set defined in (\ref{eq091}) without $Q_2$. To finish the proof of the equality in (\ref{mki7}), it remains to show that $P_1\in D(\ox,\ov)$ with $D(\ox,\ov)$ taken from (\ref{eq090}). Indeed, with $u\in\dom\sub^2\phi_\infty(\ox,\ov)$ we get from (\ref{domcod}) that $J_1\subset I_{0,1}(u)=P_1$. By $J_1\in D(\ox,\ov)$ it follows that $P_1\in D(\ox,\ov)$ and hence $(P_1,Q_1)\in{\cal A}$ with ${\cal A}$ defined in (\ref{eq092}). Employing finally Theorem~\ref{lcplw} yields $w\in\sub^2\phi_\infty(\ox,\ov)(u)$, which verifies the inclusion ``$\supset$" in (\ref{ki7}) in both cases for $u$ under consideration. To finish the proof of (ii), it suffices to combine (\ref{ki7}), (\ref{ki22}), and (\ref{ki23}). $\h$

\begin{Remark}{\bf(failure of AIQC for $\phi_\infty$).}\label{rem1} {\rm Observe that AIQC {\em fails} in the setting of Theorem~\ref{infinity}, but we still have the  {\em equality} formula for the second-order subdifferential of $\phi_\infty$. Note that the role of AICQ in the proof of Theorem~\ref{percise} is to ensure that $H_{\{Q_1,Q_2\}}\ne\emp$, $\Gamma(J_1)=J_1$, and $\Gamma(J_2)=J_2$. As Theorem~\ref{infinity} shows, this can be obtained even in the absence of AICQ.}
\end{Remark}

Next we proceed with calculating the second-order subdifferential of the 1-{\em norm function}
\begin{eqnarray}\label{1-norm}
\phi_1(x):=\|x\|_1=\disp\sum_{i=1}^n|x_i|\;\mbox{ for }\;x=(x_1,\ldots,x_n)\in\R^n.
\end{eqnarray}
Basic convex analysis tells us that $\phi_1(x)$ can be represented as the conjugate $\dd_{\B_\infty}^*(x)$ to the indicator function of the ball $\B_\infty:=\{y=(y_1,\ldots,y_n)\in\R^n|\;\|y\|_\infty\le 1\}$. Therefore
\begin{equation}\label{ki10}
u\in\sub^2\phi_1(\ox,\ov)(w)\Longleftrightarrow-w\in\sub^2\dd_{\B_\infty}(\ov,\ox)(-u)\;\mbox{ for }\;(\ox,\ov)\in\gph\partial\phi_1.
\end{equation}

The following theorem provides a precise and explicit calculation of the second-order subdifferential of $\dd_{\B_\infty}$ in terms of the initial data, which then allows us to calculate this construction for the 1-norm (\ref{1-norm}). Since in the case of $\ov\in\int\B_\infty$ we trivially have $\sub^2\dd_{\B_\infty}(\ov,\ox)(u)=\{0\}$ for any $u\in\R^n$, the emphasis below is on the boundary case where $\ov\in\bd\B_\infty$.

\begin{Theorem}{\bf (second-order subdifferential of $\dd_{\B_\infty}$).}\label{ininf}
Let $\ox\in\partial\dd_{\B_\infty}=N(\ov;\B_\infty)$ with $\ov=(\ov_1,\ldots,\ov_n)\in\bd\B_\infty$. Then we have the precise formulas:

{\bf(i)} The domain of $\partial^2\dd_{\B_\infty}(\ov,\ox)$ is calculated by
$$
\dom\sub^2\dd_{\B_\infty}(\ov,\ox)=\Big\{u=(u_1,\ldots,u_n)\in\R^n\Big|\;u_i=0\;\mbox{ for all }\;i\in I_\infty\Big\}
$$
with the index set $I_\infty:=\Big\{i\in\{1,\ldots,n\}\Big|\;\ox_i\ne 0\Big\}$.

{\bf(ii)} For any $u\in\dom\sub^2\dd_{\B_\infty}(\ov,\ox)$ the second-order subdifferential values are calculated by
$$
\sub^2\dd_{\B_\infty}(\ov,\ox)(u)=\span\Big\{e_i\Big|\;i\in E_1(u)\Big\}+\cone\Big\{\ov_i e_i\Big|\;i\in E_2(u)\Big\},
$$
where the index sets $E_1(u)$ and $E_2(u)$ are defined by
$$
E_1(u):=\Big\{i\in H(\ov)\Big|\;u_i=0\Big\}\;\mbox{ and }\;E_2(u):=\Big\{i\in H(\ov)\Big|\;\ov_i u_i>0\Big\}
$$
as some subsets of $H(\ov):=\{i\in\{1,\ldots,n\}\;\mbox{ with }\;|\ov_i|=1\}$.
\end{Theorem}
{\bf Proof.} We begin by observing the polyhedral representation of the ball in question:
\begin{eqnarray}\label{inf-ball}
\B_\infty=\Big\{y=(y_1,\ldots,y_n)\in\R^n\Big|\;\la d_i,y\ra\le 1\;\;\mbox{for}\;i=1,\ldots,2n\Big\}
\end{eqnarray}
with $d_i:=e_i\in\R^n$ and $d_{n+i}:=-e_i\in\R^n$ as $i=1,\ldots,n$. It follows from $\ov\in\bd\B_\infty$ that $I(ov)\ne\emp$ for the active constraint indexes (\ref{active}) in (\ref{inf-ball}) and that the generating vectors $\{d_i|\;i\in I(\ov)\}$ are linearly independent. By $\ox\in N(\ov;\B_\infty)$ and the normal cone representation (\ref{norc}) we find $\bar\mu_i\ge 0$ for $i\in I(\ov)$ such that $\ox=(\ox_1,\ldots,\ox_n)=\sum_{i\in I(\ov)}\bar\mu_i d_i$. This allows us to obtain by (\ref{eq05}) with $J_2:=J_+(\ox,\ov)$ the explicit expressions for the multipliers:
\begin{equation}\label{ki14}
\bar\mu_i=\left\{\begin{array}{ll}
|\ox_i|&\mbox{if }\;i\le n\\
|\ox_{i-n}|&\mbox{if }\;i>n
\end{array}\right.\quad
\mbox{for}\quad i\in J_2.
\end{equation}
Employing (\ref{ki14}) and recalling the notation $I_\infty$ from the statement in (i), we get
\begin{equation}\label{ki15}
\Big[ i\in J_2\Longrightarrow\left\{\begin{array}{ll}
i\in I_\infty&\mbox{if}\;i\le n\\
i-n\in I_\infty&\mbox{if}\;i>n
\end{array}\right.\Big],\quad
\Big[i\in I_\infty\Longrightarrow\left\{\begin{array}{ll}
i\in J_2&\mbox{if }\;\sgn(\ox_i)>0\\
i+n\in J_2&\mbox{if }\;\sgn(\ox_i)<0
\end{array}
\right.\Big].
\end{equation}

After these preparations we are ready to prove both assertions of the theorem. Since the vectors $\{d_i|\;i\in I(\ov)\}$ are linearly independent, AIQC holds and ensures by Lemma~\ref{aicq} that $\Gamma(J_2)=J_2$, which verifies (i) by using (\ref{ki15}) and Corollary~\ref{domind}. To justify (ii), we apply Theorem~\ref{percise} ensuring the equality
\begin{equation}\label{ki18}
\sub^2\dd_{\B_\infty}(\ov,\ox)(u)=\span\Big\{d_i\Big|\;i\in I_{0,2}(u)\Big\}+\cone\Big\{d_i\Big|\;i\in I_{>,2}(u)\Big\}.
\end{equation}
Combining (\ref{charcset}) and (\ref{ki14}) leads us to the relationships
\begin{equation}\label{ki16}
i\in I_{0,2}(u)\Longleftrightarrow i\in I(\ov),\quad\la d_i,u\ra=0\Longrightarrow\left\{\begin{array}{ll}
i\in E_1(u)&\mbox{if}\;i\le n\\
i-n\in E_1(u)&\mbox{if}\;i>n
\end{array}\right.,
\end{equation}
\begin{equation}\label{ki2200}
i\in E_1(u)\Longrightarrow\left\{\begin{array}{ll}
i\in I_{0,2}(u)&\mbox{if }\;{\sgn(\ov_i)>0,}\\
i+n\in I_{0,2}(u) &\mbox{if }\;{\sgn(\ov_i)<0,}\\
\end{array}\right.
\end{equation}
\begin{equation}\label{ki17}
i\in I_{>,2}(u)\Longleftrightarrow i\in I(\ov),\quad\la d_i,u\ra>0\Longrightarrow\left\{
\begin{array}{ll}
i\in E_2(u)&\mbox{if}\;i\le n\\
i-n\in E_2(u)&\mbox{if}\;i>n
\end{array}
\right.,
\end{equation}
\begin{equation}\label{ki2211}
i\in E_2(u)\Longrightarrow\left\{\begin{array}{ll}
i\in I_{>,2}(u)&\mbox{if }\;{\sgn(\ov_i)>0,}\\
i+n\in I_{>,2}(u)&\mbox{if }\;{\sgn(\ov_i)<0,}\\
\end{array}\right.
\end{equation}
which verify the formula in (ii) and thus complete the proof of the theorem. $\h$\vspace*{0.05in}

As a consequence of Theorem~\ref{ininf}, we derive now explicit formulas for calculating the domain and values of $\partial^2\phi_1$ for the 1-norm function (\ref{1-norm}). These issues have been recently addressed in \cite[Theorem~4]{e14}, where equivalent results have been derived by a different way.

\begin{Corollary} {\bf (second-order subdifferential of the 1-norm function).}\label{codnorm1} The following assertions hold in the notation Theorem~{\rm\ref{ininf}} with $(\ox,\ov)\in\gph\partial\phi_1$:

{\bf(i)} The domain of $\partial^2\phi_1(\ox,\ov)$ is calculated by
$$
\dom\sub^2\phi_1(\ox,\ov)=\Big\{w=(w_1,\ldots,w_n)\in\R^n\Big|\;w_i=0\;\mbox{ for all }\;i\in\{1,\ldots,n\}\setminus H(\ov)\Big\}.
$$

{\bf(ii)} For any $w\in\dom\sub^2\phi_1(\ox,\ov)$ we have
\begin{equation}\label{eq666}
\begin{array}{ll}
\sub^2\phi_1(\ox,\ov)(w)=\Big\{u=(u_1,\ldots,u_n)\in\R^n\Big|&u_i=0\;\;\mbox{if}\;\;i\in C_2(w)\cup I_{\infty},\\
&\ov_iu_i\le 0\;\;\mbox{if}\;\;i\in C_1(w)\setminus C_2(w)\;\Big\},
\end{array}
\end{equation}
where the index sets $C_1(w)$ and $C_2(w)$ are defined by
$$
C_1(w):=\Big\{i\in H(\ov)\Big|\;w_i\ne 0\Big\}\;\;\mbox{and }\;\;C_2(w):=\Big\{i\in H(\ov)\Big|\;w_i\ne 0,\;\ov_iw_i>0\Big\}.
$$
\end{Corollary}
{\bf Proof.} Combining (\ref{ki10}) and Theorem~\ref{ininf} tells us that $u\in\sub^2\phi_1(\ox,\ov)(w)$ amounts to
\begin{equation}\label{eq777}
\left\{\begin{array}{ll}
-u\in\Big\{u=(u_1,\ldots,u_n)\in\R^n\Big|\;u_i=0\;\mbox{ for all }\;i\in I_\infty \Big\}\\
-w\in\span\Big\{e_i\Big|\;i\in E_1(-u)\Big\}+\cone\Big\{\ov_i e_i\Big|\;i\in E_2(-u)\;\Big\}.
\end{array}\right.
\end{equation}
It follows from $E_1(-u)\cup E_2(-u)\subset H(\ov)$ that $w_i=0$ whenever $i\not\in H(\ov)$. This together with (\ref{eq777}) verifies the inclusion ``$\subset$" in (i). To justify the converse inclusion, take $w$ from the set on the right-hand side in (i). Using Theorem~\ref{ininf}, we get $w\in\sub^2\dd_{\B_\infty}(\ov,\ox)(\ou)$ for $\ou:=0\in\R^n$. This along with (\ref{ki10}) leads us the inclusion "$\supset$" in (i) and thus completes the proof of (i).

The verification of (ii) is a bit more involved. Pick $(w,u)\in\gph\sub^2\phi_1(\ox,\ov)$ and observe that $u_t=0$ for all $t\in I_{\infty}$ by (\ref{eq777}). Then for each $t\in C_2(w)$ there are two possible cases: either (a) $w_t>0$ and $\ov_t=1$, or (b) $w_t<0$ and $\ov_t=-1$. We claim that $u_t=0$ in both cases. Indeed, suppose on the contrary that $u_t\ne 0$ in case (a). It results from $w_t\ne 0$ that $t\in E_1(-u)\cup E_2(-u) $. By $u_t\ne 0$ we deduce that $t\in E_2(-u)$. Employing this together with $\ov_t=1$ yields $u_t<0$, which means that the coefficients for $e_t$ in (\ref{eq777}) are nonnegative. Thus $w_t\le 0$, a contradiction. Similar arguments lead us to a contradiction in case (b).

Next we pick $t\in C_1(w)\setminus C_2(w)$ and consider again the two possible cases: either (1) $w_t>0$ and $\ov_t=-1$, or (2) $w_t <0$ and $\ov_t=-1$. In case (1) we have $u_t\ge 0$ since $t\in E_1(-u)\cup E_2(-u)$. This shows that $\ov_t u_t\le 0$. Case (2) is treated similarly. Thus both these cases give us
$\ov_t u_t\le 0$ for all $t\in C_1(w)\setminus C_2(w)$, which justifies the inclusion ``$\subset$" in (\ref{eq666}).

To verify the opposite inclusion in (\ref{eq666}), let $u=(u_1,\ldots,u_n)$ belong to the right-hand side of (\ref{eq666}). Remembering that  $w=(w_1,\ldots,w_n)\in\dom\sub^2\phi_1(\ox,\ov)$, we need to prove that both inclusions in (\ref{eq777}) are satisfied. In fact, the first one follows immediately due to $u_i=0$ for all $i\in C_2(w)\cup I_{\infty}$. To check the second inclusion in (\ref{eq777}), fix any component $w_i\ne 0$ of $w$ for $i\in H(\ov)$ and suppose that $w_i>0$. If $\ov_i=1$ in (\ref{eq777}), then $i\in C_2(w)$ and hence $u_i=0$, which implies that $i\in E_1(-u)$. On the other hand, for $\ov_i=-1$ we get $i\in C_1(w)\setminus C_2(w)$, which yields $i\in E_1(-u)\cup E_2(-u)$. Similar arguments are applied to the case of $w_i<0$, and hence we have (\ref{eq777}). It justifies the inclusion ``$\supset$" in (\ref{eq666}) and so completes the proof of the corollary. $\h$

\end{document}